\documentclass[final,1p,times,number]{elsarticle}
\usepackage[colorlinks]{hyperref}

\usepackage{amsfonts,amssymb, mathrsfs, amsmath,   float}
\usepackage{color}
\usepackage{epsfig}

\newcommand{\SPC}{\hspace*{15pt}}
\newcommand{\qedd}{\hspace*{\fill}$\Box$\medskip}

\newtheorem{theorem}{Theorem}

\newtheorem{lemma}[theorem]{Lemma}
\newtheorem{definition}[theorem]{Definition}
\newdefinition{remark}{Remark}
\newtheorem{example}[theorem]{Example}

\newproof{pot}{\bf Proof of (*)}
\newproof{proof}{\bf Proof}

\usepackage{tikz}
\usetikzlibrary{arrows,shapes,positioning}
\usetikzlibrary{decorations.markings}
\tikzstyle arrowstyle=[scale=1]
\tikzstyle directed=[postaction={decorate,decoration={markings,mark=at position .65 with {\arrow[arrowstyle]{stealth}}}}]
\tikzstyle reverse directed=[postaction={decorate,decoration={markings,mark=at position .65 with {\arrowreversed[arrowstyle]{stealth};}}}]
\usepackage{subfigure}

\floatstyle{ruled}
\newfloat{algorithm}{t}{loa}
\floatname{algorithm}{Algorithm}
\newcommand{\Inp}[1]
  {\noindent\begin{tabular}{@{}p{1.8cm}@{}p{11.2cm}@{}}
   {\bf Input: }&#1 \end{tabular}}
\newcommand{\Outp}[1]
  {\noindent\begin{tabular}{@{}p{1.8cm}@{}p{13.2cm}@{}}
   {\bf Output: }&#1 \end{tabular}}

\def\J{\text{\rm{Jac}}}
\def\Y{{\mathbb{Y}}}
\def\A{{\mathcal A}}
\def\P{{\mathbb P}}

\def\bu{{\mathbf{u}}}
\def\bv{{\mathbf{v}}}

\def\drem{\hbox{\rm{rem}}}

\def\dzero{\mathbb{V}}

\def\sat{\hbox{\rm{sat}}}
\def\asat{\hbox{\rm{asat}}}
\def\max{\hbox{\rm{max}}}
\def\Chow{\text{Chow}}
\def\deg{\hbox{\rm{deg}}}
\def\init{\hbox{\rm{I}}}
\def\ord{\hbox{\rm{ord}}}

\def\H{\hbox{\rm{H}}}
\def\lead{\hbox{\rm{ld}}}
\def\Sep{\hbox{\rm{S}}}

\def\dim{\hbox{\rm{dim}}}

\def\mod{\hbox{\rm{mod}}}
\def\ld{{\rm ld}}

\def\coeff{\hbox{\rm{coeff}}}

\def\I{\mathcal{I}}
\def\ff{{\mathcal F}}
\def\CI{{\mathcal{I}}}
\def\ee{{\mathcal E}}

\def\trdeg{\hbox{\rm{tr.deg}}}
\def\dtrdeg{\hbox{\rm{d.tr.deg}}}
\def\Jac{\text{\rm{Jac}}}

\def\F{{\mathcal {F}}}

\def\lc{{\rm lc}}

\begin{document}
\begin{frontmatter}
\title{Computation of Differential Chow Forms for Prime Differential Ideals\tnoteref{mytitlenote}}
\tnotetext[mytitlenote]{Partially
       supported by a National Key Basic Research Project of China (2011CB302400) and  by grants from NSFC
       (60821002,11301519).}
\author{Wei Li and Yinghong Li}
\address {KLMM, Academy of Mathematics and Systems Science\\ Chinese Academy of Sciences, Beijing 100190, China}
\ead {liwei@mmrc.iss.ac.cn, liyinghong10@mails.ucas.ac.cn}
\begin{abstract}
In this paper, we propose algorithms to compute differential Chow forms for prime differential ideals which are given by their characteristic sets.
The main algorithm is based on an optimal bound for the order of a prime differential ideal in terms of its characteristic set under an arbitrary ranking,
which shows the Jacobi bound conjecture holds in this case.
Apart from the order bound, we also give a degree bound for the differential Chow form.
In addition, for prime differential ideals given by their characteristic sets under an orderly ranking, a much more simpler algorithm is given to compute its differential Chow form.
The computational complexity of both is single exponential in terms of the Jacobi number,
the maximal degree of the differential polynomials in the characteristic set and the number of variables.
\end{abstract}
\begin{keyword}
Differential Chow form \sep Jacobi bound \sep Characteristic set \sep Single exponential algorithm
\end{keyword}

\end{frontmatter}

\section{Introduction}

The Chow form, also known as the Cayley form, is a basic concept in algebraic geometry \cite{chow,hodge}
and  a powerful tool in elimination theory.
It preserves many interesting properties of the corresponding variety and also has important applications in many fields.
 For example,
Wu managed to define Chern numbers for algebraic varieties with arbitrary singularities via the Chow form \cite{wentsun1987chern};
The Chow form was also used as a tool to obtain deep results in
transcendental number theory by Nesterenko  \cite{nesterenko1977estimates}
and  Philippon  \cite{philippon1986criteres};
Brownawell made a major breakthrough in elimination theory by
developing new properties of the Chow form and proving an effective
version of the Nullstellensatz with optimal bounds
\cite{brownawell1987bounds}.
Recent study also shows that the Chow form has a close relation with sparse elimination theory \cite{gelfand,sturmfels1991sparse}.
 All these show the necessity of developing efficient algorithms to compute the Chow form.

Krick et al. showed that the Chow form of an equidimensional variety given by polynomial equations can be computed in single exponential time via an effective version of quantifier elimination in the first order theory of algebraically closed fields \cite{krick}.
Caniglia gave an algorithm which is based on linear algebra to compute the Chow form for an unmixed polynomial ideal in single exponential time and
as an application, the computational information about primary decomposition of this  ideal can be obtained \cite{Caniglia1990compute}.
Jeronimo et al. gave a bounded probabilistic algorithm whose complexity
is polynomial in the size and the geometric degree of the input polynomial
equation system \cite{Jeronimo2004computational}.

Differential algebraic geometry founded by Ritt and Kolchin aims at studying differential equations in a similar way that polynomial equations are studied in
algebraic geometry \cite{ritt1966differential, kolchin1973differential}.
Recently, we generalized the algebraic Chow form to its differential analog and  the theory of differential Chow forms in both affine and
projective differential algebraic geometry was developed
\cite{gao2010intersection,li2011differential}. It has been shown that most of the basic
properties of the algebraic Chow form can be extended to its
differential counterpart \cite{gao2010intersection}.
Next, it is quite natural to explore further problems related to differential Chow form in both algorithmic and applied aspects.
By its definition, we know the differential Chow form can be
computed by the characteristic set method.
However, it is hard to estimate the computational complexity if using this method.
Recall that the differential Chow form preserves the main properties of the differential ideal,
so these properties should be considered in order to realize efficient algorithms for computing the differential Chow form.
Bearing this principle in mind, in this paper, we focus on devising efficient algorithms to compute differential Chow forms for prime differential ideals.

In general, there is no algorithm to test whether a given ideal is prime or not.
However, for most applications, prime differential ideals are often given by their characteristic sets with respect to some ranking.
Thus, the main problem we consider in this paper is as follows.
Given a prime differential ideal represented by its characteristic set under an arbitrary ranking,  devise an algorithm to compute its differential Chow form and estimate the computational complexity in the worst case.
Although this can be realized by means of algorithms on transforming characteristic sets from one ranking to another \cite{boulier2010computing, sadik2000bound, golubitsky2009algebraic} as explained in \cite{gao2010intersection},
these algorithms either lack complexity analysis or are so general as not efficient enough.
We will propose an algorithm to compute the differential Chow form in single exponential time which requires only linear algebraic computations in the base field of the ideal.

For a prime differential ideal $\mathcal{I}=\sat(\mathcal{A})$ in the differential polynomial ring $\ff\{y_1,\ldots,y_n\}$, where $\mathcal{A}$ is a characteristic set of $\mathcal{I}$ w.r.t. some ranking $\mathscr{R}$, the dimension of $\CI$ can be read off from $\mathcal{A}$ which is just equal to $n-|\mathcal{A}|$.
Intersecting $\mathcal{I}$ with $n-|\mathcal{A}|+1$ generic differential hyperplanes,
 by definition,  the differential Chow form is just a differential polynomial  in the coefficients of these hyperplanes with minimal order and also of minimal degree under this order among all polynomials contained in the intersection ideal.
Naturally, we first need to give  bounds for the order and the degree of the differential Chow form.
First, by \cite[Theorem 4.11]{gao2010intersection}, the order of the differential Chow form is just equal to the order of the corresponding prime differential ideal.
So it is equivalent to give an order bound for the prime differential ideal in terms of its characteristic set.
Here, two cases should be considered according to whether $\mathscr{R}$ is an orderly ranking or not.
If $\mathcal{A}$  is a characteristic set of $\CI$  under some orderly ranking,
then the precise order of $\mathcal{I}$ is just equal to $\ord(\mathcal{A})$.
But when $\mathscr{R}$  is an arbitrary ranking, the problem becomes much more complicated.
In \cite{golubitsky2009algebraic}, Golubitsky et al. obtained an order bound by proposing $\ord(\CI)\leq |\mathcal{A}|\cdot\max_{A\in\mathcal{A}}\ord(A)$.
This bound is non-optimal, and they conjectured that $\ord(\CI)\leq\sum_{i=1}^{|\mathcal{A}|}m_{i}$, where $m_i=\max_{A\in\mathcal{A}}\ord(A,y_i)$ and  $m_{1}\geq m_{2}\geq\cdots\geq m_{n}$, without giving a proof.
In this paper, combining the result of Kondrativa et al. on Jacobi's bound for systems of  independent differential polynomials, we prove that the order of $\CI$
is bounded by the Jacobi number of $\mathcal{A}$, which is a better bound than that in the above conjecture as shown in Example~\ref{bound}.

We also give a Be\'{z}out-type degree bound for the prime differential ideal in terms of the degrees of the differential polynomials in its characteristic set.
Then based on the order and degree bounds, we give algorithms to compute differential Chow forms.
The algorithms require only linear algebraic computations in the base field of the ideals and
the computational complexities in the worst case are single exponential in terms of Jacobi numbers.

The rest of the paper is organized as follows. In section 2, some basic notation and preliminary results about differential algebra will be given.
Section 3 contributes to give an algorithm to compute differential Chow forms for prime differential ideals represented by characteristic sets w.r.t. orderly rankings.
In section 4, for more general prime differential ideals given by characteristic sets under arbitrary rankings,
we give an algorithm to compute differential Chow forms.
Finally, we conclude this paper and propose open problems for further research in section 5.

\section{Preliminaries}
In this section, some basic notations and preliminary results in differential algebra will be given. For more details about differential algebra, please refer to \cite{buium1999differential,kondratieva1998differential,kolchin1973differential,ritt1966differential,sit2002ritt}.

\subsection{Differential polynomial algebra}
Let $\mathcal {F}$ be a fixed ordinary differential field of characteristic zero with a derivation operator $\delta$.
For ease of notations,  we use primes and exponents $(i)$ to denote derivatives  under $\delta$,
and for each $a\in \mathcal {F}$, denote $a^{[n]}=\{a,a^{(1)},\ldots,a^{(n)}\}$ and $a^{[\infty]}=\{a^{(i)}|i\geq0\}$.
Throughout this paper, unless otherwise indicated, $\delta$ is kept fixed during any discussion.
A typical example of differential field is $\mathbb{Q}(t)$ which is the field of rational functions in the variable $t$ with $\delta=\frac{d}{dt}$.

Let $\mathcal {G}$ be a differential extension field of $\mathcal {F}$ and $S$ a subset of $\mathcal {G}$.
 We denote respectively by $\mathcal {F}[S]$, $\mathcal {F}(S)$, $\mathcal {F}\{S\}$, and $\mathcal {F}\langle S\rangle$ the smallest subring, the smallest subfield, the smallest differential subring, and the smallest differential subfield of $\mathcal {G}$ containing $\mathcal {F}$ and $S$.
And $\mathcal{G}$ is said to be finitely generated over $\mathcal{F}$ if there exists a finite subset $S\subset\mathcal{G}$ such that $\mathcal {G}=\mathcal {F}\langle S\rangle$.

Let $\Theta$ be the free communicative semigroup with unit (written multiplicatively) generated by $\delta.$
A subset $\Sigma$ of a differential extension field $\mathcal {G}$ of $\mathcal {F}$ is said to be differentially  dependent  over $\mathcal {F}$ if the set $(\theta \alpha)_{\theta \in \Theta, \alpha \in \Sigma}$ is algebraically dependent over $\mathcal {F}$, and otherwise, it is said to be differentially  independent  over $\mathcal {F}$, or to be a family of differential  indeterminates  over $\mathcal {F}$ (abbr. differential $\mathcal {F}$-indeterminates). In the case $\Sigma$ consists of only one element $\alpha$, we say that $\alpha$ is differentially algebraic or differentially transcendental over $\mathcal {F}$ respectively.
A maximal subset $\Omega$ of $\mathcal {G}$ which is differentially  independent over $\mathcal {F}$ is said to be a differential transcendence basis of $\mathcal {G}$ over $\mathcal {F}$.
We use $\dtrdeg\,\mathcal {G}/\mathcal {F}$ to denote the differential transcendence  degree of $\mathcal {G}$ over $\mathcal {F}$, which is the cardinal number of $\Omega$.

Suppose $\mathcal{G}_1$ and $\mathcal{G}_2$ are two differential extension fields of $\F$.
A homomorphism (reps. isomorphism)  $\phi$ from $\mathcal{G}_1$ to $\mathcal{G}_2$ is called {\em a differential homomorphism (reps. isomorphism) over $\F$} if $\phi$ commutes with $\delta$ and leaves each element of $\F$ invariant.

A differential extension field $\mathcal {E}$ of $\mathcal {F}$ is called a  universal  differential extension field,
if for any finitely generated differential extension field $\mathcal {F}_1\subset\mathcal {E}$ of $\mathcal {F}$ and any finitely generated differential extension field $\mathcal {F}_2$ of $\mathcal {F}_1$ not necessarily contained in $\mathcal {E}$,  there exists a differential extension field $\mathcal{F}_3\subset\mathcal {E}$  of
$\mathcal {F}_1$  such that $\mathcal{F}_3$ is differentially isomorphic to $\mathcal {F}_2$ over $\mathcal {F}_1$. Such a universal differential extension field of $\mathcal {F}$ always exists \cite{kolchin1973differential}.

Now suppose $\mathcal {E}$ is a universal differential extension field of $\mathcal {F}$, and $\mathbb{Y}=\{y_1,\ldots,y_n\}$ is a set of differential indeterminates over $\mathcal {E}$.
 For any $y\in\Y$, denote $\delta^ky$ by $y^{(k)}.$ The elements of $\mathcal
{F}\{\Y\}=\mathcal {F}[y_j^{(k)}\,|\,j=1,\ldots,n;k\in \mathbb{N}]$
are called {\em differential polynomials} over $\ff$ in $\Y$,
 and $\mathcal {F}\{\Y\}$ itself is called the {\em differential polynomial ring } over $\ff$ in $\Y$.
 A differential polynomial ideal $\mathcal {I}$ in $\mathcal {F}\{\mathbb{Y}\}$ is an  algebraic ideal which is closed under derivation, i.e. $\delta(\mathcal {I})\subseteq \mathcal {I}$.
  And a prime differential ideal is a differential ideal which is prime as an algebraic polynomial ideal.
  For convenience, a prime differential ideal is assumed not to be the unit ideal in this paper.

By a differential affine space, we mean any one of the sets $\mathcal {E}^n(n\in \mathbb{N})$.
Let $\Sigma$ be a subset of differential polynomials in $\mathcal {F}\{\mathbb{Y}\}$.
 A point $\eta\in\mathcal{E}^n$ is called a differential zero of $\Sigma$ if $f(\eta)=0$ for any $f\in \Sigma$.
 The set of all differential zeros of $\Sigma$ is denoted by $\mathbb{V}(\Sigma)$, which is called a differential variety defined over $\mathcal {F}$.
 A point $\eta\in \mathbb{V}(\mathcal {I})$ is called a  generic point  of a prime differential ideal $\mathcal {I} \subseteq \mathcal {F}\{\mathbb{Y}\}$
 if for any differential polynomial $f\in \mathcal {F}\{\mathbb{Y}\}$ we have $f(\eta)=0 \Leftrightarrow f\in \mathcal {I}$. It is well known that:

\begin{lemma}
A non-unit differential ideal is prime if and only if it has a generic point.
\end{lemma}

We conclude this section by giving the definition and some basic properties of K{\"a}hler differentials which will be used in this paper.
(See \cite{Eisenbud1995commutative} for more details on K{\"a}hler differentials in the purely algebraic case and \cite{johnson1969kahler} for  K{\"a}hler differentials in differential algebra.)

\begin{definition}
Let $R$ be a field and  $S$  an $R$-algebra. The  module of {\em  K{\"a}hler  differentials}  of $S$ over $R$, written $\Omega_{S/R}$, is the $S$-module generated by the set $\{d(f): f\in S\}$ subject to the relations
$$d(bb')=bd(b')+b'd(b)$$
$$d(ab+a'b')=ad(b)+a'd(b')$$
for all $a,a'\in R$, and $b,b'\in S.$
\end{definition}

\begin{theorem}\cite{johnson1969kahler} \label{th-algindep}
If $R$ is a field of characteristic zero and $S$ is a field extension of $R$, then the elements $\eta_1,\ldots,\eta_r$ of $S$ are algebraically independent over $R$ if and only if $d(\eta_1),\ldots,d(\eta_r)$ are linearly independent over $S$.
\end{theorem}

\begin{lemma}\cite{johnson1969kahler}\label{ka}
If $\ff$ is a differential field and $S$ is a differential algebra over $\ff$ with derivation operator $\delta$, then $\Omega_{S/\ff}$ has a canonical structure of differential module over $S$ such that  for each $b\in S$,
$$\delta d(b)=d\delta(b).$$
\end{lemma}

\subsection{Characteristic sets of a differential polynomial ideal}
Let $f$ be a differential polynomial in $\mathcal {F}\{\mathbb{Y}\}$.
The order of $f$ w.r.t. $y_i$ is the greatest number $k$ such that $y_i^{(k)}$ appears effectively in $f$, denoted by $\ord(f,y_i)$.
If $y_i$ does not appear in $f$,  set $\ord(f,y_i)=-\infty$.
The $order$ of $f$ is defined to be $\max_i\{\ord(f,y_i)\}$,  denoted by $\ord(f)$.

A $ranking$ $\mathscr{R}$ is a total order over $\Theta(\mathbb{Y})$ if satisfying
1) $\delta\alpha>\alpha$ for all $\alpha\in \Theta(\mathbb{Y})$
and
2) $\alpha_1>\alpha_2 \Rightarrow \delta\alpha_1>\delta\alpha_2$ for all $\alpha_1, \alpha_2\in \Theta(\mathbb{Y})$.
Below are two important kinds of ranking:
\begin{enumerate}[1.]
\item $Elimination$ $ranking$: $y_i>y_j\Rightarrow \delta^k y_i>\delta^l y_j$ for any $k,l \geq 0$.
\item $Orderly$ $ranking$: $k>l \Rightarrow \delta^k y_i>\delta^l y_j$ for any $i,j\in\{1,\ldots,n\}$.
\end{enumerate}
Let $f$ be a differential polynomial in $\mathcal {F}\{\mathbb{Y}\}$ endowed with a ranking $\mathscr{R}$.
The {\em leader} of $f$ is the greatest derivative w.r.t. $\mathscr{R}$ which appears effectively in $f$, denoted by $\lead(f)$.
Regarding $f$ as a univariate polynomial in $\ld(f)$, its leading coefficient is called the {\em initial} of $f$, denoted by $\init_f$,
and the partial derivative of $f$ w.r.t. $\ld(f)$ is called the  {\em separant} of $f$, denoted by $\Sep_f$.
For any two differential polynomials $f$, $g$ in $\mathcal {F}\{\mathbb{Y}\}$, $f$ is said to be of  lower rank than $g$, denoted by $f<g$, if 1) $\ld(f)<\ld(g)$, or 2) $\ld(f)=\ld(g)$ and $\deg(f,\ld(f))<\deg(g,\ld(f))$. And $f$ is said to be reduced w.r.t. $g$ if no proper derivatives of $\ld(g)$ appear in $f$ and $\deg(f,\ld(g))<\deg(g,\ld(g))$. Let $\mathcal {A}$ be a set of differential polynomials. Then $\mathcal {A}$ is said to be an auto-reduced set  if each differential polynomial in $\mathcal {A}$ is reduced w.r.t any other element of $\mathcal {A}$. Every auto-reduced set is finite \cite{ritt1966differential}.

Let $\mathcal {A}$ be an auto-reduced set. We denote $\H_\mathcal {A}$ to be the set of all initials and separants of $\mathcal {A}$ and $\H_\mathcal {A}^\infty$  the minimal multiplicative set containing $\H_\mathcal {A}$. The  saturation differential ideal  of $\mathcal {A}$ is defined to be
$$\sat(\mathcal {A})=[\mathcal {A}]:\H_\mathcal {A}^\infty=\{f\in \mathcal {F}\{\mathbb{Y}\}\mid \exists h\in \H_\mathcal {A}^\infty, \text{s.t.} hf\in [\mathcal {A}]\}.$$ The algebraic saturation ideal of $\mathcal{A}$ is defined to be $\asat(\mathcal {A})=(\mathcal {A}):\init_\mathcal {A}^\infty$, where
$\init_\mathcal{A}^\infty$ is the multiplicative set generated by the initials of polynomials in $\mathcal{A}$.
 We use capital calligraphic letters such as $\mathcal{A}, \mathcal{B}, \ldots$ to denote auto-reduced sets 
 and use notation $\mathcal{A}=A_1,\ldots, A_t$ to specify the list of the elements of $\mathcal{A}$ arranged in order of increasing rank.

 Let $\mathcal {A}=A_1,\ldots,A_t$ be an auto-reduced set with $\init_i$ and $\Sep_i$ being the initial and separant of $A_i$,
 and $f$ an arbitrary differential polynomial. 
 Then there exists an algorithm, called Ritt-Kolchin algorithm of differential reduction \cite{sit2002ritt}, 
 which reduces $f$ w.r.t. $\mathcal{A}$ to a differential polynomial $r$ that is reduced w.r.t. $\mathcal{A}$, satisfying  that
$$\prod_{i=1}^t\Sep_i^{d_i}\init_i^{e_i}\cdot f\equiv r, \mod\, [\mathcal {A}],$$
where $d_i$ and $e_i\,(i=1,\ldots,t)$ are nonnegative integers.
We call this $r$ the {\em differential remainder} of $f$ w.r.t. $\mathcal {A}$,
denoted by $\drem(f,\mathcal{A})$.
Throughout the paper, the differential remainder w.r.t. an auto-reduced set is  always assumed to be computed by performing the reduction algorithm described in \cite[Section 6]{sit2002ritt}.

An auto-reduced set $\mathcal {C}$ contained in a differential polynomial set $\mathcal {S}$ is said to be a $characteristic$ $set$ of $\mathcal {S}$ if $\mathcal {S}$ does not contain any nonzero element reduced w.r.t $\mathcal {C}$. A characteristic set $\mathcal {C}$ of a differential ideal $\mathcal {J}$ reduces all elements of $\mathcal {J}$ to zero. Furthermore, if $\mathcal {J}$ is prime, then $\mathcal {J}=\sat(\mathcal {C})$.

\begin{definition}
For an auto-reduced set $\mathcal{A}= A_1, \ldots, A_t$ with $\ld(A_i)=y_{c_i}^{(o_i)}$, the order of $\mathcal{A}$ is defined as $\ord(\mathcal {A})=\sum_{i=1}^t o_i$, and the set $\mathbb{Y}\setminus\{y_{c_1},\ldots,y_{c_t}\}$ is called the {\em parametric set} of $\mathcal{A}$.
\end{definition}

Finally, we recall the definition of  differential dimension and order for a prime differential ideal $\CI$, which are
closely related to characteristic sets of $\CI$.

 Let $\mathcal {I}$ be a prime differential ideal in
$\ff\{\Y\}$ and $\xi=(\xi_{1},\ldots,\xi_{n})$ a generic point of
$\mathcal {I}$.
The {\em differential dimension} of $\mathcal {I}$ or $\dzero(\mathcal {I})$
is defined as the differential transcendence degree of  the
differential extension field $\mathcal {F}\langle
\xi_{1},\ldots,\xi_{n}\rangle$ over $\mathcal {F}$, that is,
$\dim(\mathcal {I})=\dtrdeg\, \mathcal {F}\langle
\xi_{1},\ldots,\xi_{n}\rangle/\mathcal {F}$.
By \cite{hubert2000factorization}, the differential dimension of $\mathcal{I}$ is equal to the cardinal number of the
parametric set of its characteristic set w.r.t. any ranking.

%

\begin{definition}\cite{kolchin1964notion}
Let  $\mathcal {I}$ be a prime differential  ideal of $\mathcal
{F}\{\Y\}$ with a generic point $\eta=(\eta_{1},\ldots,\eta_{n})$.
Then there exists a unique numerical polynomial $\omega_{\mathcal
{I}}(t)$  such that $\omega_{\mathcal {I}}(t)=\trdeg\,\mathcal
{F}(\eta_{i}^{(j)}:i=1,\ldots,n;j\leq t)/\mathcal {F} $ for all
sufficiently large $t \in \mathbb{N}$. $\omega_{\mathcal {I}}(t)$ is
called the {\em differential dimension polynomial} of $\mathcal {I}$.
\end{definition}

\begin{theorem}\cite[Theorem 13]{sadik2000bound}\label{ord-equ}
Let $\mathcal {I}$ be a prime differential ideal of dimension $d$.
Then the differential dimension polynomial of $\mathcal{I}$ is of  the form $\omega_{\mathcal {I}}(t)=d(t+1)+h$.
The number $h$ is defined to be the order of $\mathcal {I}$, denoted by $\ord(\mathcal {I})$.
And if $\mathcal{A}$ is a characteristic set of $\mathcal {I}$ under any orderly ranking, then $\ord(\mathcal {I})=\ord(\mathcal{A})$.
\end{theorem}

\subsection{Chow form for a prime differential ideal}
In this section, we recall the definition of the differential Chow form and some of its basic properties. For more details about differential Chow form, please refer to \cite{gao2010intersection}.

A generic differential hyperplane is the zero set  of a differential polynomial $u_0+u_1y_1+\cdots+u_ny_n$  contained in $\ee^n$ where the $u_i\in\mathcal{E}$ are differentially independent over $\mathcal {F}$.
Let  $\mathcal {I}\subseteq \mathcal {F}\{\mathbb{Y}\}$ be a prime differential ideal of dimension $d$ and $$\mathbb{P}_i=u_{i0}+u_{i1}y_1+\cdots +u_{in}y_n\, (i=0,\ldots,d)$$ be $d+1$ generic differential hyperplanes.
For each $i$, denote $\bu_i=(u_{i0},u_{i1},\ldots,u_{in})$,  and $\bu=\bigcup_{i=0}^d \bu_i\setminus \{u_{i0}\}$.
 Let $[\mathcal {I},\mathbb{P}_0,\ldots,\mathbb{P}_d]$ be the ideal generated by $\CI$ and the $\P_i$ in $\ff\{\Y,\bu_0,\ldots,\bu_d\}.$
Then by \cite{gao2010intersection},
\begin{lemma}\label{chow-1}
$[\mathcal {I},\mathbb{P}_0,\ldots,\mathbb{P}_d]\cap \mathcal {F}\{\bu_0,\ldots,\bu_d\}$ is a prime differential ideal of codimension one.
\end{lemma}

By the  theory of characteristic sets,  there exists an irreducible differential polynomial $F(\bu_{0},$ $\ldots,\bu_{d})\in\mathcal {F}\{ \bu_{0},\ldots,\bu_{d}\}$ such that $\{F\}$ is a characteristic set of $[\mathcal {I},\mathbb{P}_0,\ldots,\mathbb{P}_d]\cap \mathcal {F}\{\bu_0,\ldots,\bu_d\}$ under
any ranking. That is, $$[\mathcal {I},\mathbb{P}_0,\ldots,\mathbb{P}_d]\cap \mathcal {F}\{\bu_0,\ldots,\bu_d\}=\sat(F).$$
This $F$ is defined to  be the {\em differential Chow form} of $\mathcal {I}$.

Differential Chow forms can uniquely characterize their corresponding differential ideals.
The following theorem gives some basic properties of differential Chow forms which will be used later.

\begin{theorem} \label{th-chowproperty}
Let  $\CI\subset\mathcal{F}\{\Y\}$ be a prime differential ideal of dimension $d$ and order $h$ 
with $F(\bu_0,\bu_1,\ldots,\bu_d)$ its  differential Chow form.
Then the following assertions hold.
\begin{enumerate}[1.]
\item $\ord(F)=h$. And for each $u_{ij}$ appearing effectively in $F$, we have $\ord(F,u_{ij})=h$. In particular, $\ord(F,u_{i0})=h$.
\item $F(\bu_0,\ldots\bu_d)$ is differentially homogeneous  of the same degree in each $\bu_i\,(i=0,\ldots,d)$.
Namely, there exists a nonnegative integer $r$ such that for each $i$ and  a differential indeterminate $\lambda$,
$F(\bu_0,\ldots, \lambda\bu_i,\ldots,\bu_d)=r\cdot F(\bu_0,\ldots, \bu_i,\ldots,\bu_d)$.
\item $\A=F,\frac{\partial F}{\partial u_{00}^{(h)}}y_1-\frac{\partial F}{\partial u_{01}^{(h)}},\ldots,\frac{\partial F}{\partial u_{00}^{(h)}}y_n-\frac{\partial F}{\partial u_{0n}^{(h)}}$
is a characteristic set
\footnote{Here $\A$ is a  triangular set but may not be an ascending chain.
Note that the differential remainder of  $\frac{\partial F}{\partial u_{00}^{(h)}}$ w.r.t. $F$ is not zero,
so $\A$ can also serve as a characteristic set.}
of $[\mathcal {I},\P_0,\ldots,\P_d]\subset\mathcal {F}\{\bu_{0},\ldots,\bu_{d},$ $\Y\}$
w.r.t. the elimination ranking $u_{d0}\prec\cdots\prec u_{00}\prec y_1\prec\cdots\prec y_n$.
\item Suppose $F_{\rho\tau}$ is obtained from $F$ by interchanging $\bu_\rho$ and $\bu_\tau$ in $F$. Then $F_{\rho\tau}$ and $F$ differ at most by a sign.
\end{enumerate}
\end{theorem}

\section{Computation of differential Chow forms  for prime differential ideals represented by characteristic sets under orderly rankings}
In this section, we will give an algorithm to compute the differential Chow form for a prime differential ideal represented by its characteristic set w.r.t. an orderly ranking based on linear algebraic techniques and then estimate its computational complexity.

Given a prime differential ideal $\sat(\mathcal{A})$ with $\mathcal{A}$ its characteristic set under an orderly ranking,
by Theorems~\ref{ord-equ} and \ref{th-chowproperty},  the order of its differential Chow form  is  equal to $\ord(\mathcal{A})$.
In order to estimate the computational complexity, degree bounds are also needed.
So before giving the algorithm, we first give the degree bound for the differential Chow form.

\subsection{Degree bound of the differential Chow form in terms of characteristic set under an orderly ranking}

 In this section, we will give the degree bound for the differential Chow form of a prime differential ideal.
 Before giving the main result, we firstly need some lemmas.

\begin{lemma}\label{lm-elimination} \cite{Heintz1983definability,li2011sparse}
Let $\I$ be a prime ideal in $\mathcal {F}[\Y]$ and $\I_r=\I\cap \ff[x_1,\ldots,x_r]$ for any $1\leq r\leq n$.
Then $\deg(\I_r)\leq \deg(\I)$.
\end{lemma}

\begin{lemma}\cite[Proposition 1]{lazard1983grobner} \label{le-deg-pol}
Let $F_1,\ldots,F_m \in \mathcal {F}[\Y]$ be  polynomials
generating an ideal $\I$ of dimension $r$. Suppose
$\deg(F_1)\geq\cdots\geq\deg(F_m)$ and let
$D:=\prod_{i=1}^{n-r}\deg(F_i)$. Then $\deg(\I)\leq D$.
\end{lemma}

\begin{lemma}\label{ele1}
Let $\mathcal {I}\subseteq \mathcal {F}\{\mathbb{Y}\}$ be a prime differential ideal of dimension $d$ with $\{A_1,\ldots,A_{n-d}\}$  as a characteristic set of $\mathcal {I}$ w.r.t. an orderly ranking and $e_i=\ord(A_i)$, $h=\sum_{i=1}^{n-d}e_i$.
Suppose $F$ is the differential Chow form of $\mathcal {I}$. Then  $$(F)=(A_1^{[h-e_1]},\ldots,A_{n-d}^{[h-e_{n-d}]},\mathbb{P}_0^{[h]},\ldots,\mathbb{P}_d^{[h]}, \H x_0-1)\cap \mathcal {F}[\bu_0^{[h]},\ldots,\bu_d^{[h]}],$$
where $\H=\prod_{i=1}^{n-d}\init_{A_i}\Sep_{A_i}$ and $x_0$ is a new indeterminant.
\end{lemma}

\proof
First, we claim that $\mathcal {I}\cap \mathcal {F}[\mathbb{Y}^{[h]}]=(A_1^{[h-e_1]},\ldots,A_{n-d}^{[h-e_{n-d}]},\H x_0-1)\cap \mathcal {F}[\mathbb{Y}^{[h]}].$
 Let $J=(A_1^{[h-e_1]},\ldots,$ $A_{n-d}^{[h-e_{n-d}]}, \H x_0-1)\bigcap \mathcal {F}[\mathbb{Y}^{[h]}]=\asat(A_1^{[h-e_1]},\ldots,A_{n-d}^{[h-e_{n-d}]})\bigcap \mathcal {F}[\mathbb{Y}^{[h]}]$.
For any $f\in \mathcal {I}\cap \mathcal {F}[\mathbb{Y}^{[h]}]$, there exists an integer $l\geq 0$ such that $H^lf=\sum_{i=1}^{n-d}\sum_{k_i=0}^{h-e_i}g_{k_i}A_i^{(k_i)}=
[( \H x_0-1+1)/{x_0}]^lf$, where $g_{k_i}\in \mathcal {F}[\mathbb{Y}^{[h]}]$.
So $f\in(A_1^{[h-e_1]},\ldots,A_{n-d}^{[h-e_{n-d}]}, \H x_0-1)$, and consequently $f\in J$.
On the other hand, for any $f\in J$, we have $f=\sum_{i=1}^{n-d}\sum_{k_i=0}^{h-e_i}g_{k_i}A_i^{(k_i)}+g( \H x_0-1)$, here $g_{k_i}, g\in \mathcal {F}[\mathbb{Y}^{[h]},x_0]$. Thus if we substitute $x_0=1/H$ into this equality, we get $f\in \mathcal {I}\cap \mathcal {F}[\mathbb{Y}^{[h]}]$. Hence $\mathcal {I}\cap \mathcal {F}[\mathbb{Y}^{[h]}]=J$.

Thus, we have
\begin{eqnarray}
(F)&=&\big(\mathcal {I}\cap \mathcal {F}[\mathbb{Y}^{[h]}],\mathbb{P}_0^{[h]},\ldots,\mathbb{P}_d^{[h]}\big)\cap \mathcal {F}[\bu_0^{[h]},\ldots,\bu_d^{[h]}]
\nonumber \\
&=&\Big((A_1^{[h-e_1]},\ldots,A_{n-d}^{[h-e_{n-d}]}, \H x_0-1)\cap \mathcal {F}[\mathbb{Y}^{[h]}],\mathbb{P}_0^{[h]},\ldots,\mathbb{P}_d^{[h]}\Big)\bigcap \mathcal {F}[\bu_0^{[h]},\ldots,\bu_d^{[h]}]
\nonumber \\
&\subseteq&(A_1^{[h-e_1]},\ldots,A_{n-d}^{[h-e_{n-d}]},\mathbb{P}_0^{[h]},\ldots,\mathbb{P}_d^{[h]}, \H x_0-1)\cap \mathcal {F}[\bu_0^{[h]},\ldots,\bu_d^{[h]}]
\nonumber \\
&\subseteq&[A_1,\ldots,A_{n-d},\mathbb{P}_0,\ldots,\mathbb{P}_d, \H x_0-1]\cap \mathcal {F}[\bu_0^{[h]},\ldots,\bu_d^{[h]}]
\nonumber \\
&=&[\CI,\P_0,\ldots,\P_d]\cap\mathcal {F}[\bu_0^{[h]},\ldots,\bu_d^{[h]}]\nonumber\\
&=& (F)
\nonumber
\end{eqnarray}
Thus, $(F)=(A_1^{[h-e_1]},\ldots,A_{n-d}^{[h-e_{n-d}]},\mathbb{P}_0^{[h]},\ldots,\mathbb{P}_d^{[h]}, \H x_0-1)\cap \mathcal {F}[\bu_0^{[h]},\ldots,\bu_d^{[h]}]$.
\qedd

The following theorem gives the degree bound of the differential Chow form.

\begin{theorem}\label{ord-deg-bound}
Let $\mathcal {I}\subseteq \mathcal {F}\{\mathbb{Y}\}$ be a prime differential ideal of dimension $d$ and $\{A_1,\ldots,A_{n-d}\}$  its characteristic set w.r.t. an orderly ranking. Suppose $F$ is the differential Chow form of $\mathcal {I}$.
Let $e_i=\ord(A_i)$, $h =\sum_{i=1}^{n-d}e_i$ and
$\deg(A_i)=m_i$, then $$\deg(F)\leq 2^{(h+1)(d+1)}\prod_{i=1}^{n-d}m_i^{h-e_i+1}(2\sum_{i=1}^{n-d}(m_i-1)+1).$$
In particular, let $D=\max\{2,m_i\}$, then
$\deg(F)<(n-d+1)D^{(nh+n+3)}.$
\end{theorem}

\proof
Set $\mathcal {J}=(A_1^{[h-e_1]},\ldots,A_{n-d}^{[h-e_{n-d}]},\mathbb{P}_0^{[h]},\ldots,\mathbb{P}_d^{[h]}, \H x_0-1)$, where $\H=\prod_{i=1}^{n-d}\init_{A_i}\Sep_{A_i}$. By Lemma~\ref{le-deg-pol}, we have $\deg(\mathcal {J})\leq (\prod_{i=1}^{n-d}m_i^{h-e_i+1})2^{(h+1)(d+1)}\big(2\sum_{i=1}^{n-d}(m_i-1)+1\big)$.
And by Lemmas~\ref{lm-elimination} and \ref{ele1}, we have $\deg(F)=\deg\big(\mathcal {J}\bigcap \mathcal {F}[\bu_0^{[h]},\ldots,\bu_d^{[h]}]\big)\leq \deg(\mathcal {J})$. Thus, \begin{eqnarray}
\deg(F)&\leq&2^{(h+1)(d+1)}\prod_{i=1}^{n-d}m_i^{h-e_i+1}(2\sum_{i=1}^{n-d}(m_i-1)+1))\nonumber\\
&<& (n-d+1)D^{(nh+n+3)}\nonumber
\end{eqnarray}
\qedd

\subsection{Complexity of differential reductions}

In this section, we estimate the complexity of performing differential reductions, which will be very useful when analyzing the computational complexity of differential Chow form. Before doing so, we first recall some results about algebraic reductions.

\begin{lemma}\cite[Lemma 3.3.3]{Szanto}\label{alg-1}
Let $f\in \ff[x_{1},\ldots,x_{l}]\,(l\geq n)$ and $g\in \ff[x_{1},\ldots,x_{n}]$ with $m=\deg(g,x_n)>0$.
Regarding $f$ and $g$ as univariate polynomials in $x_{n}$,
 suppose we have computed $r\in \ff[x_1,\ldots,x_l]$ with $\deg_{x_n}(r)<\deg_{x_n}(g),$ such that there exits $q\in \ff[x_1,\ldots,x_l]$ satisfying that $$(\lc(g,{x_n}))^{k+1}f=qg+r, $$
 where $k=\deg_{x_n}(f)-m$ and $\lc(g,x_n)$ is the coefficient of $g$ in $x_n^m$.
 Then for $j<n$, $$\deg_{x_j}(q)\leq (k+1)\deg_{x_j}(g)+\deg_{x_j}(f),$$ $$\deg_{x_j}(r)\leq (k+2)\deg_{x_j}(g)+\deg_{x_j}(f),$$ and for $j>n$, $$\deg_{x_j}(q), \deg_{x_j}(r)\leq \deg_{x_j}(f).$$
\end{lemma}

The above $r$ is called the {\em algebraic pseudo-remainder} of $f$ w.r.t. $g$.
Based on the above lemma, we can now analyze the computational complexity of reducing a polynomial w.r.t.
an autoreduced set.

\begin{lemma}\label{alg-2}
Let $\mathcal {A}=\{A_1,\ldots,A_p\}$ be an auto-reduced set in $\ff[x_1,\ldots,x_n]$ w.r.t. any fixed ranking $\mathcal {R}$ with $\ld(A_i)=y_{i}\,(i=1,\ldots,p)$. Set $m=\max_i\{\deg(A_{i})\}$.
Then for any $f\in \mathcal {F}[x_1,\ldots,x_n]$, $\deg(f)=D$,
the algebraic pseudo-remainder $r$ of $f$ w.r.t. $\mathcal {A}$ can be computed with at most  $(p+1)^{2.376}[(m+1)^p(D+1)]^{2.376n}$
 $\mathcal {F}$-arithmetic operations and the degree of $r$ is bounded by $(m+1)^p(D+1)-1$.
\end{lemma}

\proof Let $f_p=f, f_{p-1},\ldots,f_0=r$ be the pseudo-remainder sequence of $f$ w.r.t. $\mathcal{A}$  satisfying the following equations:
$$(\init_{A_i})^{l_i}f_{i}=q_iA_i+f_{i-1},$$ where
$$l_i=\deg_{y_{i}}(f_i)-\deg_{y_{i}}(A_i)+1.$$
Then, by Lemma~\ref{alg-1}, $\deg(f_i)$ and $\deg(q_i)$ satisfy the following relations:
$$\deg(f_p)=\deg(f)=D,$$
$$\deg(f_{i-1})\leq (m+1)\deg(f_i)+m,$$
$$\deg(q_i)\leq (m+1)\deg(f_i).$$
So for all $0\leq i\leq p$, $\deg(f_i)\leq (m+1)^{p-i}(D+1)-1$, and for all $1\leq j\leq p$, $\deg(q_j)\leq (m+1)^{p-j+1}(D+1)-m-1$.
It follows that $\deg(r)\leq (m+1)^p(D+1)-1.$

From above, we have the following expression
$$\prod_{i=1}^p\init_{A_i}^{l_i}f=\prod_{i=1}^{p-1}\init_{A_i}^{l_i}q_pA_p+\prod_{i=1}^{p-2}\init_{A_i}^{l_i}q_{p-1}A_{p-1}+\cdots+\init_{A_1}^{l_1}q_2A_2+q_1A_1+r.$$
Regarding this expression as a polynomial equation in the $x_i$ and  collecting the coefficients of distinct monomial terms,
we can get a system of linear polynomial equations  over $\ff$  in coefficients of $q_i$ and $r$ whose degree bounds are given above.
Thus,  $r$ can be computed by solving this linear equation system consisting of at most
$w_1\leq{(m+1)^p(D+1)+n\choose n}$ equations in  $w_2=\sum_{i=1}^{p}{\deg(q_i)+n\choose n}+{\deg(r)+n\choose n}$ variables.
To solve it, we need at most $(\max\{w_1,w_2\})^\omega$ $\mathcal {F}$-arithmetic operations, where $\omega$ is the matrix multiplication exponent and the currently best known $\omega$ is 2.376. Here, $w_1\leq [(m+1)^p(D+1)]^n$ and $w_2\leq\sum_{i=1}^p\binom{(D+1)(m+1)^{p-i+1}-1+n}{n}+\binom{(D+1)(m+1)^{p}-1+n}{n}]\leq (p+1)[(m+1)^p(D+1)]^n$.
So the computional complexity is
\begin{eqnarray}
\big(\max\{w_1,w_2\}\big)^\omega &\leq & (p+1)^{2.376}[(m+1)^p(D+1)]^{2.376n}.
\nonumber
\end{eqnarray}\qedd

Let $f$ and $g$ be two differential polynomials in $\mathcal {F}\{y_1,\ldots,y_n\}$.
Suppose $\ld(g)=y_\alpha^{(o)}$.
Since the differential remainder of $f$ w.r.t. $g$ is just equal to the algebraic remainder of $f$ w.r.t. $\{g, g^{(1)}, \ldots, $ $g^{(l)}\}$ where $l=\ord(f,y_{\alpha})-o$,
the computational complexity of differential reductions can be estimated by performing a series of algebraic reductions.

\begin{theorem}\label{Diff-red}
Let $\mathcal{A}=\{A_1,\ldots,A_p\}$ be a differential auto-reduced set  in $\mathcal {F}\{y_1,\ldots,y_n\}$
under some fixed ranking $\mathscr {R}$ and $f\in \mathcal {F}\{y_1,\ldots,y_n\}$.
Set $h=\ord(f)$, $D=\deg(f)$, $e=\max_i\{\ord(A_i)\}$ and $m=\max_i\{\deg(A_{i})\}$.
Then the differential remainder of $f$ w.r.t. $\mathcal{A}$ can be computed with at most $O\big(p(h+1)\big)(m+1)^{O\big(np(h+1)(e+h+1)\big)}(D+1)^{O(n(e+h+1))}\big)$
$\mathcal {F}$-arithmetic operations and its degree is bounded by $(m+1)^{p(h+1)}(D+1)-1$.
\
\end{theorem}
\proof
Set \\
\[\mathcal{A}_{h}=\left\{\begin{array}{l}A_{1},A_{1}^{(1)}\ldots,A_{1}^{(h)}\\ \cdots \\ A_{p},A_{p}^{(1)},\ldots,A_{p}^{(h)} \end{array}
\right.\]
Set $l=p(h+1)$, $N=n(e+h+1)$. Rewrite $\mathcal{A}_h$ to be an ascending triangular set w.r.t. the ordering induced by $\mathscr {R}$
and denote it by $\mathcal{B}=\{B_1, B_2, \ldots, B_l\}$.
Then to compute the differential remainder of $f$ with respect to $\mathcal{A}$, it suffices to compute the algebraic pseudo remainder of $f$ with respect to $\mathcal{B}$.
By Lemma~\ref{alg-2},  the differential remainder of $f$ w.r.t. $\mathcal{A}$ can be computed with at most $(l+1)^{2.376}[(m+1)^l(D+1)]^{2.376N}= [p(h+1)+1]^{2.376}[(m+1)^{p(h+1)}(D+1)]^{2.376n(e+h+1)}\leq O\big(p(h+1)\big)(m+1)^{O\big(np(h+1)(e+h+1)\big)}(D+1)^{O(n(e+h+1))}\big)$ $\mathcal {F}$-arithmetic operations.
\qedd

\subsection{An algorithm to compute the differential Chow form}

Let $\mathcal {I}=\sat(\mathcal{A})\subset\ff\{\Y\}$ be a prime differential ideal of dimension $d$  and $\mathcal{A}=\{A_1,\ldots,A_{n-d}\}$ is a characteristic set of $\mathcal{I}$ w.r.t. an orderly ranking $\mathscr{R}$.
Let $\mathbb{P}_{i}=u_{i0}+u_{i1}y_1+\cdots+u_{in}y_n\,(i=0,\ldots,d)$ be $d+1$ generic differential hyperplanes.
Set $\bu_i=\{u_{i0},\ldots,u_{in}\}$ and $\bu=\cup_{i=0}^d\bu_{i}\backslash \{u_{i0}\}$.
Let $\mathscr{R}_1$ be the elimination ranking with $\Y< \bu< u_{00}<\cdots<u_{d0}$ and $\mathscr{R}_1|_{\Y}=\mathscr{R}$.
By \cite{gao2010intersection},
$\{\mathcal {A},\mathbb{P}_{0},\ldots,\mathbb{P}_{d}\}$ is a characteristic set of the prime differential ideal $[\mathcal {I},\mathbb{P}_{0},\ldots,\mathbb{P}_{d}]\subset\mathcal {F}\{\bu_0,\ldots,\bu_d,\mathbb{Y}\}$ w.r.t. the ranking $\mathscr{R}_1$.
Moreover, if $F$ is the differential Chow form of $\mathcal {I}$,  then $[\mathcal {I},\mathbb{P}_0,\ldots,\mathbb{P}_d]\cap \mathcal {F}\{\bu_0,\ldots,\bu_d\}=\sat(F)$ and $\ord(F)=\ord(\mathcal{A})$.
Therefore, if $F_0$ is a homogeneous differential polynomial of the smallest degree among all polynomials in $\mathcal {F}[\bu_i^{(k)}:i=0,\ldots,d;0\leq k\leq \ord(\mathcal{A})]$ whose differential remainder w.r.t. $\{\mathcal {A},\mathbb{P}_0,\ldots,\mathbb{P}_d\}$ under $\mathscr {R}_1$ is zero, then $F_0$ must be the differential Chow form of $\mathcal {I}$.

With the above idea, now we give an algorithm to compute the differential Chow form of $\mathcal {I}$.
With the fixed order $h=\ord(\mathcal{A})$, the algorithm works adaptively by searching $F$ from degree $t=1$.
 If we cannot find an $F$ with such a degree, then we repeat the procedure with degree $t+1$.
 Theorem \ref{ord-deg-bound} guarantees the termination of the algorithm with  $t\leq2^{(h+1)(d+1)}\prod_{i=1}^{n-d}\deg(A_i)^{h-\ord(A_i)+1}\cdot$
 $(2\sum_{i=1}^{n-d}(\deg(A_i)-1)+1)$. In this way, we need only to handle problems with the real size and need not go to the upper bound in most cases.

\begin{algorithm}[htbp]\label{alg-diffChowform1}
  \caption{\bf --- DChowForm-1($\mathcal {A}$)} \smallskip
  \Inp{A characteristic set $\mathcal {A}=\{A_{1},\ldots,A_{n-d}\}$ of a prime differential ideal $\mathcal {I} $
under    an orderly ranking $\mathscr{R}$}\\
  \Outp{The differential Chow form $F(\bu_0,\ldots,\bu_{d})$ of $\mathcal {I}$.}\medskip
  \noindent

  1. For $i=0,\ldots,d$, let
  $\P_i=u_{i0}+u_{i1}y_1+\cdots+u_{in}y_n$ and
  $\bu_i=(u_{i0},\ldots,u_{in})$.\\
  2. Set $h=\ord(\mathcal {A})$. Set $\bv=\cup_{i=0}^{d}\bu_i^{[h]}$.\\
  3. Set $F=0$ and $t=1$.\\
  4. While $F=0$ do\\
   \SPC 4.1.  Set $F_0$ to be a homogenous GPol of degree $t$ in $\bv$.\\
   \SPC 4.2.    Set $\textbf{c}=\coeff(F_0,\bv)$. \\
   \SPC 4.3. Substitute $u_{i0}^{(k)}=-(u_{i1}y_1+\cdots+u_{in}y_n)^{(k)}\,(i=0,\ldots,d; k\geq0)$ into $F_{0}$ to get $F_1$.\\
  \SPC 4.4. Compute $F_2=\drem(F_1,\mathcal {A})$ under ranking $\mathscr{R}$.\\
   \SPC 4.5.  Set $\mathcal{P}=\coeff(F_2,\Theta(\Y)\cup\bv)$. Note $\mathcal{P}$ is a set of linear homogenous polynomials in $\textbf{c}$.\\
   \SPC 4.6.  Solve the linear equation system $\mathcal{P}=0$.\\
   \SPC 4.7. If $\textbf{c}$ has a non-zero solution, then substitute it into $F_{0}$ to get $F$ and return $F$;\\
      \SPC \SPC \, else $F=0$.\\
   \SPC  4.8.  $t:=t+1$.\medskip
  \noindent

  /*/\; Pol and GPol stand for algebraic polynomial
  and generic  algebraic polynomial.\smallskip
   \noindent

   /*/\; $\coeff(F,V)$ returns the set of coefficients of $F$ as an algebraic  polynomial in
$V$.
\smallskip
\noindent

/*/\; $\drem(f, \mathcal{B})$ returns the differential remainder of $f$ w.r.t. an auto-reduced set $\mathcal{B}$.
\end{algorithm}

\begin{theorem}\label{com-1}
Let $\mathcal {I}=\sat(\mathcal{A})$ be a prime differential ideal of dimension $d$ and $\mathcal{A}=\{A_1,\ldots,A_{n-d}\}$ is a differential characteristic set of $\mathcal{I}$ under some orderly ranking.  Set $e_i=\ord(A_i)$, $h=\sum_{i}e_i$, $e=\max_i\{e_i\}$, $m_i=\deg(A_i)$ and $m=\max\{m_i\}$.
Algorithm  DChowForm-1 computes the differential Chow form $F(\bu_0,\ldots,\bu_d)$ of $\mathcal {I}$ with
at most $$O\big([n(m+1)^{O((h+1)(2n+d+1))}]^{O(n(e+dh+2h+d+1))}\big)$$ $\ff$-arithmetic operations.
\end{theorem}

\proof
The algorithm finds a nonzero differential polynomial  $F\in \mathcal {F}\{\bu_0,\ldots,\bu_d\}$ of the smallest degree satisfying that $\ord(f)=h=\ord(\mathcal{A})$
and  the differential remainder of $f$ w.r.t. $\mathcal {A},\mathbb{P}_0,\ldots,\mathbb{P}_d$ under $\mathscr{R}_1$ is zero.
The existence of such an $F$ is obvious since $[\mathcal{I},\P_0,\ldots,\P_d]\cap\ff[\bu_0^{[h]},\ldots,\bu_d^{[h]}]=(\Chow(\mathcal{I}))$,
where $\Chow(\mathcal{I})$ is the differential Chow form of $\mathcal{I}$, and this $F$ must be the differential Chow form of $\mathcal{I}$.

We estimate the computational complexity of the algorithm below.
In each loop of step 4, the complexity of the algorithm is clearly dominated by step 4.4 and step 4.6.
 In step 4.4, we need to compute the differential remainder $F_2$ of $F_1$ w.r.t. the characteristic set $\mathcal {A}$.
 By Theorem~\ref{Diff-red}, $F_2$ can be computed with at most $O\big((n-d)(h+1)(m+1)^{O(n(n-d)(h+1)(e+h+1))}(2t+1)^{O(n(e+h+1))}\big)$
 $\ff$-arithmetic operations, and the degree of $F_2$ is bounded by $(m+1)^{(n-d)(h+1)}(2t+1)-1$. In step 4.6, we need to solve the linear equation system $\mathcal {P}=0$ in $\textbf{c}$. It is easy to see that $\mid\textbf{c}\mid=\binom{(d+1)(h+1)(n+1)+t-1}{t}$, then $\mathcal {P}=0$ is a linear equation system with $W_{1,t}=\binom{(d+1)(h+1)(n+1)+t-1}{t}$ variables and $W_{2,t}=\binom{d_{F_2}+n(e+h+1)+(d+1)(h+1)(n+1)}{d_{F_2}}$ equations.
 To solve it, we need at most $\max\{W_{1,t},W_{2,t}\}^\omega$ $\mathcal {F}$-arithmetic operations, where $\omega$ is the matrix multiplication exponent and currently, the best known $\omega$ is 2.376.

Suppose $T$ is the degree bound of the differential Chow form.
The iteration in Step 4 may loop from 1 to $T$ in the worst case.
Thus, in terms of $T$, the differential Chow form can be computed with at most
\begin{eqnarray}
& &\sum_{t=1}^{T}\Big\{O\big((n-d)(h+1)(m+1)^{O(n(n-d)(h+1)(e+h+1))}(2t+1)^{O(n(e+h+1))}\big)+ W_{2,t}^{2.376} \Big\}
\nonumber\\
&\leq&T \{O\big((n-d)(h+1)(m+1)^{O(n(n-d)(h+1)(e+h+1))}(2T+1)^{O(n(e+h+1))}\big)+
\nonumber \\
& &O\big([(m+1)^{(n-d)(h+1)}(2T+1)-1]^{O(2.376(n(e+h+1)+(d+1)(h+1)(n+1)))}\big)\}
\nonumber\\
&=&O\big([(2T+1)(m+1)^{(n-d)(h+1)}]^{O(n(e+dh+2h+d+1))}\big).
\nonumber
\end{eqnarray}
$\ff$-arithematic operations. Here, to derive  the above inequalities, we always assume that $(m+1)^{(n-d)(h+1)}(2T+1)>n(e+dh+2h+d+1)$.
Hence, the theorem follows by simply replacing $T$ by the degree bound for $F$ given in Theorem~\ref{ord-deg-bound}.
\qedd

We use the following example to illustrate the above algorithm.
\begin{example}
Let $n=1$ and $\mathcal {A}=\{y'-4y\}$. Clealy, $d=\dim(\sat(\mathcal{A}))=0$.
We use this simple example to illustrate Algorithm 1.
Let $\mathbb{P}_0=u_{00}+u_{01}y$, and $\bu_0=(u_{00}, u_{01})$.
In step 2, $h=\ord(\mathcal {A})=1$,
$\bv=\bu_0^{[1]}=(u_{00},u_{01},u_{00}',u_{01}')$.
We first execute steps 4.1 to 4.7 for $t=1$.
Set $F_0=c_{01}u_{00}+c_{02}u_{01}+c_{03}u_{00}'+c_{04}u_{01}'$, $\textbf{c}=(c_{01},c_{02},c_{03},c_{04})$.
 In step 4.3., we get  $F_1=-c_{01}u_{01}y+c_{02}u_{01}-c_{03}u_{01}'y-c_{03}u_{01}y'+c_{04}u_{01}'$.
 And in Step 4.4.,  $F_2=-(c_{01}+4c_{03})u_{01}y+c_{02}u_{02}-c_{03}u_{01}'y+c_{04}u_{01}'$.
 Then $\mathcal {P}=\{c_{01}+4c_{03}, c_{02}, c_{03}, c_{04}\}$.
 Hence $\mathcal {P}=0$ has a unique solution $\textbf{c}=(0,0,0,0)$.
 In Step 4.8, $t=2$. Next we execute steps 4.1 to 4.7 for $t=2$.
 Set $F_0=c_{01}u_{00}^2+c_{02}u_{00}u_{01}+c_{03}u_{00}u_{00}'+c_{04}u_{00}u_{01}'+c_{05}u_{01}2+c_{06}u_{01}u_{00}'
+c_{07}u_{01}u_{01}'+c_{08}u_{00}'^2+c_{09}u_{00}'u_{01}'+c_{10}u_{01}'^2$, $\textbf{c}=(c_{01}, \ldots, c_{10})$. In step 4.3 and step 4.4, we get  $F_1=c_{01}u_{01}^2y^2-c_{02}u_{01}^2y+c_{03}u_{01}u_{01}'y^2+c_{03}u_{01}^2yy'-c_{04}u_{01}u_{01}'y+c_{05}u_{01}^2
-c_{06}u_{01}u_{01}'y-c_{06}u_{01}^2y'+c_{07}u_{01}u_{01}'+c_{08}(u_{01}'y+u_{01}y')^2-c_{09}(u_{01}'^2y)-c_{09}u_{01}u_{01}'y'
+c_{10}u_{01}'^2$, and $F_2=(c_{01}+4c_{03}+16c_{08})u_{01}^2y^2+(c_{03}+8c_{08})u_{01}u_{01}'y^2+c_{08}u_{01}'2y^2-(c_{02}+4c_{06})u_{01}^2y
-(c_{04}+c_{06}+4c_{09})u_{01}u_{01}'y-c_{09}u_{01}'2y+c_{05}u_{01}^2+c_{07}u_{01}u_{01}'+c_{10}u_{01}'^2$ respectively.
So $\mathcal {P}=0$ consists of equations $\{c_{01}+4c_{03}+16c_{08}=0, c_{03}+8c_{08}=0, c_{05}=c_{07}=c_{08}=c_{09}=c_{10}=0, c_{02}+4c_{06}=0, c_{04}+c_{06}+4c_{09}=0\}$. Hence $\textbf{c}=(0,4q,0,q,0,-q,0,0,0,0)$ where $q\in\mathbb{Q}$.
Substitute $\textbf{c}$ into $F_0$, then we get $F=4u_{00}u_{01}+u_{00}u_{01}'-u_{01}u_{00}'.$
Therefore, this algorithm returns $F=4u_{00}u_{01}+u_{00}u_{01}'-u_{01}u_{00}'$, which is exactly the differential Chow form of $\CI=\sat(\mathcal{A})$.
\end{example}

\section{Computation of differential Chow forms for differential ideals represented by characteristic sets under arbitrary differential rankings}
In this section, we will consider the general case. Namely, for prime differential ideals represented by characteristic sets under arbitrary rankings, we give algorithms to compute the differential Chow form.

\subsection{Order bound of the differential Chow form}
%
%

In the case that a prime differential ideal is represented by a characteristic set under an orderly ranking, the order of the Chow form is equal to the order of this characteristic set. However, this may not be true for arbitrary rankings. In this section, we will give an upper bound for the order of the differential Chow form of a prime differential ideal in terms of the Jacobi number of its characteristic set w.r.t any fixed ranking.

Let $\mathcal{S}=\{f_1,\ldots,f_n\}$ be $n$ differential polynomials in $\Y$.
Let $e_{ij} = \ord(f_i,y_j)$ be the order of $f_i$ in $y_j$ if $y_j$
occurs effectively in $f_i$ and $e_{ij} = -\infty$ otherwise.
Then the {\em Jacobi bound}, or the {\em Jacobi number}, of $\mathcal{S}$,
denoted as $\Jac(\mathcal{S})$, is the maximum of the summations of
all the diagonals of $E=(e_{ij})$. Or equivalently,
$$\Jac(E) = \max_\sigma\sum_{i=1}^n e_{i\sigma(i)},$$
where $\sigma$ is a permutation of $\{1,\ldots,n\}$.
The {\em Jacobi's Problem} conjectures that the order of every zero
dimensional component of $\mathcal{S}$ is bounded by the Jacobi number of
$\mathcal{S}$ \cite{ritt1935jacobi}. This conjecture is closely related to the differential dimension conjecture \cite{cohn1983order}
which was also proposed by Ritt \cite{ritt1966differential}.
Both the two well-known conjectures in differential algebra still remain open and they are proved only in some special cases,
for instance, $n=1$ or linear polynomial systems.

In the latest two decades, many differential algebraists  work on the Jacobi's order bound conjecture and proposed
other order bounds for prime differential ideals in
terms  of either generators or characteristic sets under arbitrary rankings \cite{kondrat2008jacobi, golubitsky2009algebraic}.
Let $\mathcal{I}=\sat(\mathcal{A})$ be a prime differential ideal with $\mathcal{A}$ as a characteristic set under an arbitrary ranking.
In \cite{golubitsky2009algebraic}, Golubitsky et al.  showed $\ord(\mathcal{I})\leq |\mathcal{A}|\cdot\max\{\ord(C): C\in\mathcal{A}\}$.
Moreover, since this bound is likely to be not optimal, they  proposed another better order bound in terms of Ritt number \cite{ritt1966differential}.
Namely, let $o_i=\max_{C\in\mathcal{A}}\{\ord(C,y_i)\}$ for $i=1,\ldots,n$, where set $\ord(C,y_i)=0$ if $y_{i}$ does not occur in $C$.
Suppose $o_{k_1}\geq o_{k_2}\geq\cdots\geq o_{k_n}$ is arranged in non-increasing order, then
they conjectured $\ord(\mathcal{I})\leq\sum_{i=1}^{|\mathcal{A}|}o_{k_i}$ without giving a proof.
Clearly, the Jacobi bound is optimal to this conjectured bound.

As a main result of this section, we will prove that Jacobi's order bound conjecture holds for prime differential ideals
specified by characteristic sets.  Firstly, we recall some results from \cite{kondrat2008jacobi} for later use.
\begin{definition}
Suppose that $P$ is a prime ideal of a commutative ring $B$ and $M$ is a $B$-module. A set $H\subseteq M$ is called independent over $P$ if $\{h+PM\mid h\in H\}$ is a system of elements of $M/MP$ linearly independent over the quotient ring $B/P$.
\end{definition}

\begin{definition} \label{def-independent}
Let $\I$ be a prime differential ideal of the differential polynomial ring $\F\{\Y\}$.
The set $\{f_1,\ldots,f_k\}\subset\F\{\Y\}$ is called independent over $\mathcal {I}$ if the set $\{df_{j}^{(i)}\mid 1\leq j\leq k,i\geq 0\}\subset\Omega_{\F\{\Y\}/\mathcal {F}}$ is independent over $\mathcal {I}$.
Here, $(\Omega_{\F\{\Y\}/\mathcal {F}},d)$ is the module of  K{\"a}hler differentials of $\F\{\Y\}$ over $\mathcal {F}$ as defined in Lemma~\ref{ka}.
\end{definition}

\begin{lemma}\label{ord1}\cite[Theorem 3]{kondrat2008jacobi}
Let $\mathcal {I}$ be a prime differential ideal in $\mathcal {F}\{y_{1},\ldots,y_{n}\}$.
 Suppose $f_{1},\ldots,f_{n}\in\I$.
 If $f_{1},\ldots,f_{n}$ are independent over $\mathcal {I}$,
 then $\ord(\mathcal {I})\leq \J(f_{1},\ldots,f_{n})$.
\end{lemma}

The following lemma is crucial to prove our result about the order bound of the differential Chow form.

\begin{lemma}\label{ord2}
Let $\mathcal {I}\subset \mathcal {F}\{\mathbb{Y}\}$ be a prime differential ideal of dimension $d$, and $\mathcal{A}=\{A_{1},\ldots,A_{n-d}\}$ be the characteristic set under any fixed ranking $\mathscr{R}$.
Let $L_{i}=u_{i0}+u_{i1}y_1+\cdots+u_{in}y_n(i=1,\ldots,d)$ be $d$ independent generic differential hyperplanes with coefficient vector $\bu_i=(u_{i0},\ldots,u_{in})$, and ${\mathcal {J}}=[\mathcal {I},L_{1},\ldots,L_{d}]_{\F\langle\bu_1,\ldots,\bu_d\rangle\{\Y\}}$.
Then $A_{1},\ldots,A_{n-d},L_{1},\ldots,L_{d}$ are independent over ${\mathcal {J}}$.
\end{lemma}
\proof
For convenience, suppose  $\lead(A_i)=y_{d+i}^{(o_i)}\,(i=1,\ldots,n-d)$ with $A_i< A_j\,(i<j)$ and the parametric set of $\mathcal{A}$ is $\{y_{1},\ldots,y_d\}$.
Denote $\F_d=\F\langle\bu_1,\ldots,\bu_d\rangle.$
By \cite[Theorem 3.6]{gao2010intersection}, ${\mathcal {J}}=[\mathcal {I},L_{1},\ldots,L_{d}]_{\F_d\{\Y\}}$ is a prime differential ideal.

By Definition~\ref{def-independent}, we need to show that the set $\{d(A_{j}^{(i)})+{\mathcal {J}}\Omega_{R/\mathcal {F}}, d(L_{k}^{(i)})+{\mathcal {J}}\Omega_{R/\mathcal {F}}:1\leq j\leq n-d, 1\leq k\leq d, i\geq 0\}\subset\Omega_{\F_d\{\Y\}/\F_d}$ is linearly independent over ${\mathcal {J}}$.
By Theorem~\ref{th-algindep}, it is easy to derive that $d(y_j^{{(i)}})+{\mathcal {J}}\Omega_{R/\mathcal {F}}\,(j=1,\ldots,n; i\geq0)$ are linearly independent over ${\mathcal {J}}$. So it suffices to prove that for each $k\geq0$,
the Jacobi submatrix of $\mathcal{S}_k=\{A_{1}^{[k]},\ldots,A_{n-d}^{[k]},L_{1}^{[k]},\ldots,L_{d}^{[k]}\}$ w.r.t. $\Y$ and its derivatives has full row rank module ${\mathcal {J}}$.
 Let $T$ be the $({n(k+1))\times (n(k+1)})$ submatrix of this Jacobi matrix with columns indexed by monomials $y_{d+1}^{(o_1)},\ldots,y_{n}^{(o_{n-d})},\ldots,y_{d+1}^{(o_1+k)},\ldots,$ $y_n^{(o_{n-d}+k)},y_1,\ldots,y_d,\ldots,y_1^{(k)},\ldots,y_d^{(k)}$. Then $T$ can be written in the following block form:
$$T=\left(
\begin{array}{cccccccc}
M_1 & \textbf{0}  & \cdots & \textbf{0}  & *  & * & \cdots & * \\
 \textbf{0}   & M_1         & \cdots & \textbf{0}  & *  & * & \cdots & * \\
\vdots & \vdots   & \ddots & \vdots&\vdots & \vdots & \vdots & \vdots  \\
  * &   *         &   \cdots &M_1 & * & * & \cdots & *\\
*&  * & * & * & M_2 & \textbf{0}  & \cdots & \textbf{0} \\
*&  * & * & * & \textbf{0} & M_2  & \cdots & \textbf{0} \\
\vdots & \vdots   & \vdots & \vdots&\vdots & \vdots & \ddots & \vdots  \\
*&  * & * & * & \textbf{0} & \textbf{0}   & \cdots & M_2\\
\end{array}
\right),$$
where $M_1=\left(\begin{array}{cccc}\Sep_{A_{1}} & 0& \cdots & 0\\
0 & \Sep_{A_{2}}& \cdots & 0 \\
\vdots & \vdots      & \ddots & \vdots\\
 * &   *         &   \cdots & \Sep_{A_{n-d}}\end{array}\right)$ and
  $M_2=\left(\begin{array}{cccc} u_{11} & u_{12}&\cdots & u_{1d}\\
                                  u_{21} & u_{22}&\cdots & u_{2d}\\
                               \cdots & \cdots   & \cdots & \cdots\\
                           u_{d1} & u_{d2}&\cdots & u_{dd}\\ \end{array}\right)$.

We claim that $T$ has full row rank module $\mathcal {J}$. It suffices to show that $\det(T)\notin{\mathcal {J}}$.
We first claim that for each $f \in {\mathcal {J}}\bigcap \mathcal {F}\{\widetilde{\bu},\mathbb{Y}\}$, where $\widetilde{\bu}=\{u_{ij}| 1\leq i \leq d, 1\leq j\leq n\}$,
 if we rewrite $f$ as a differential polynomial in $\widetilde{\bu}$ with coefficients in $\F\{\Y\}$,
 that is, $f=\sum_{\phi}\phi(u)g_{\phi}(\mathbb{Y})$, then $g_{\phi}(\mathbb{Y})\in \mathcal {I}$ for all $\phi$.
Indeed, let $\mathcal{J}_0=[\mathcal {I},L_{1},\ldots,L_{d}]_{\F\{\Y,\bu_1,\ldots,\bu_d\}}$ and $\xi=(\xi_{1},\ldots,\xi_{n})$ be a generic point of $\mathcal {I}$ free from $\F\langle\bu_1,\ldots,\bu_d\rangle$.
Let $\zeta=(\xi,-\sum_{i=1}^nu_{1i}\xi_{i},u_{11},\ldots,u_{1n},\ldots,-\sum_{i=1}^nu_{di}\xi_{i},u_{d1},\ldots,u_{dn})$.
It is easy to show that $\zeta$ is a generic point of $\mathcal{J}_0$
and ${\mathcal {J}}\bigcap \mathcal {F}\{\widetilde{\bu},\mathbb{Y}\}\subset\mathcal{J}_0$.
So $f(\zeta)=0$, and consequently, for each $\phi$, $g_\phi(\xi)=0$, which implies that  $g_{\phi} \in \mathcal {I}$.
Rewrite $\det(T)$ as a differential polynomial in $\widetilde{\bu}$ and suppose $\det(T)=\sum_{\phi}\phi(u)g_{\phi}(\mathbb{Y})$.
By the claim,  it remains to show that there exists a differential monomial $\phi^*(\widetilde{\bu})$ such that $g_{\phi^*}(\mathbb{Y})\notin \mathcal {I}$.
By the structure of $T$, we can take $\phi^*=(\prod_{i=1}^du_{ii})^k$ and $g_{\phi^*}=(\prod_{i=1}^{n-d}s_{A_{i}})^k$, which completes the proof.
\qedd

\begin{theorem}\label{Jac-bou}
Let $\mathcal {I}\subset \mathcal {F}\{\mathbb{Y}\}$ be a prime differential ideal of dimension $d$, and $\mathcal {A}=\{A_{1},\ldots,$ $A_{n-d}\}$ the characteristic set of $\mathcal {I}$ under any fixed ranking $\mathscr{R}$. Then $\rm{ord}(\mathcal {I})\leq \J(\mathcal {A})$.
\end{theorem}
\proof By \cite[Theorem 3.6]{gao2010intersection}, ${\mathcal {J}}=[\mathcal {I},L_{1},\ldots,L_{d}]_{\F\langle\bu_1,\ldots,\bu_d\rangle\{\Y\}}$ is a prime differential ideal with $\ord(\mathcal{J})=\ord(\mathcal{I})$.
And from Lemma~\ref{ord1} and Lemma~\ref{ord2},  $\ord(\mathcal{J})\leq  \J(\mathcal {A},L_{1},\ldots,L_{d})=\J(\mathcal {A})$.
Thus, it follows that $\rm{ord}(\mathcal {I})\leq \J(\mathcal {A})$.
\qedd

The following simple example shows that Jacobi bound is optimal to the conjectured order bound in \cite{golubitsky2009algebraic}.
\begin{example}\label{bound}
Let $\mathcal {I}=\sat(y_2y_3+1, y_1^{(n)}y_2^{(n)}+y_3)\subset\ff\{y_1,y_2,y_3\}$ be a prime differential ideal with $\{y_2y_3+1, y_1^{(n)}y_2^{(n)}+y_3\}$ a characteristic set of $\mathcal {I}$ w.r.t. the elimination ranking $y_3<y_2<y_1$.
By Theorem \ref{Jac-bou}, $\ord(\mathcal{I})\leq n$. While by the conjecture in \cite{golubitsky2009algebraic},
$\ord(\mathcal {I})\leq 2n$.
\end{example}

\subsection{Degree bound of the differential Chow form in terms of characteristic sets under arbitrary rankings}

In this section, based on the order bound given in the preceding section, we will give a degree bound for the differential Chow form of a prime differential ideal in terms of its characteristic set under an arbitrary ranking. The method used here is similar to that in Section 3.1.

\begin{lemma}\label{ele10}
Let $\mathcal {I}$ be a prime differential ideal in $\mathcal {F}\{\mathbb{Y}\}$ of dimension $d$ and $\mathcal {A}=\{A_1,\ldots,$ $A_{n-d}\}$ its characteristic set w.r.t. an arbitrary  ranking.
Suppose $F$ is the differential Chow form of $\mathcal {I}$ and $\ord(\mathcal {I})=h$. Then $$(F)=(A_1^{[h]},\ldots,A_{n-d}^{[h]},\mathbb{P}_0^{[h]},\ldots,\mathbb{P}_d^{[h]}, \H x_0-1)\cap \mathcal {F}[\bu_0^{[h]},\ldots,\bu_d^{[h]}],$$ where $H=\prod_{i=1}^{n-d}\init_{A_i}\Sep_{A_i}$ and $x_0$ is a new indeterminant.
\end{lemma}

\proof
Note that for each $f\in\mathcal {F}[\mathbb{Y}^{[h]}]$, the differential remainder of $f$ w.r.t. $\mathcal{A}$ can be obtained by computing the algebraic remainder of $f$
w.r.t. $\{A_1^{[h]},\ldots,A_{n-d}^{[h]}\}$.
 So similarly to the proof of Lemma~\ref{ele1},  it is easy to show that $\mathcal {I}\cap \mathcal {F}[\mathbb{Y}^{[h]}]=(A_1^{[h]},\ldots,A_{n-d}^{[h]}, \H x_0-1)\cap \mathcal {F}[\mathbb{Y}^{[h]}]$.
Then
\begin{eqnarray}
& &(F)\nonumber\\
&=&\big(\mathcal {I}\cap \mathcal {F}[\mathbb{Y}^{[h]}],\mathbb{P}_0^{[h]},\ldots,\mathbb{P}_d^{[h]}\big)\cap \mathcal {F}[\bu_0^{[h]},\ldots,\bu_d^{[h]}]
\nonumber \\
&=&\big((A_1^{[h]},\ldots,A_{n-d}^{[h]}, \H x_0-1)\cap \mathcal {F}[\mathbb{Y}^{[h]}],\mathbb{P}_0^{[h]},\ldots,\mathbb{P}_d^{[h]}\big)\cap \mathcal {F}[\bu_0^{[h]},\ldots,\bu_d^{[h]}]
\nonumber \\
&\subseteq&(A_1^{[h]},\ldots,A_{n-d}^{[h]},\mathbb{P}_0^{[h]},\ldots,\mathbb{P}_d^{[h]}, \H x_0-1)\cap \mathcal {F}[\bu_0^{[h]},\ldots,\bu_d^{[h]}]
\nonumber \\
&\subseteq&[A_1,\ldots,A_{n-d},\mathbb{P}_0,\ldots,\mathbb{P}_d]\cap \mathcal {F}[\bu_0^{[h]},\ldots,\bu_d^{[h]}]
\nonumber \\
&=& (F).
\nonumber
\end{eqnarray}
Thus, the lemma is proved.
\qedd

\begin{theorem}\label{any-deg-bound}
Let $\mathcal {I}$ be a prime differential ideal in $\mathcal {F}\{\mathbb{Y}\}$ of dimension $d$ and order $h$ with $\mathcal {A}=\{A_1,\ldots,$ $A_{n-d}\}$ its characteristic set under  an arbitrary  ranking. Let $F$ be the differential Chow form of $\mathcal{I}$. Set $\deg(A_i)=m_i$.  Then $$\deg(F)\leqslant (\prod_{i=1}^{n-d}m_i^{h+1})2^{(h+1)(d+1)}\big(2\sum_{i=1}^{n-d}(m_i-1)+1\big).$$ In particular, let $D=\max\{m_i,2\}$, then $\deg(F)<D^{(\Jac(\mathcal {A})+1)(n+1)}(2(D-1)(n-d)+1)$.
\end{theorem}

\proof
Let $\mathcal {J}=(A_1^{[h]},\ldots,A_{n-d}^{[h]},\mathbb{P}_0^{[h]},\ldots,\mathbb{P}_d^{[h]}, \H x_0-1)\subset\ff[\Y^{[h]},\bu_0^{[h]},\ldots,\bu_d^{[h]},x_0]$,
then by Lemma~\ref{le-deg-pol}, we have $\deg(\mathcal {J})\leq \prod_{i=1}^{n-d}m_i^{h+1}2^{(h+1)(d+1)}\big(2\sum_{i=1}^{n-d}(m_i-1)+1\big)$. From Lemma~\ref{ele10} and Lemma~\ref{lm-elimination}, we get $\deg\big((F)\big)=\deg\big(\mathcal {J}\bigcap \mathcal {F}(\bu_0^{[h]},\ldots,\bu_d^{[h]})\big)\leq \deg(\mathcal {J})$.
Thus $\deg(F)\leqslant 2^{(h+1)(d+1)}\prod_{i=1}^{n-d}m_i^{h+1}\big(2\sum_{i=1}^{n-d}(m_i-1)+1\big)$. By Theorem~\ref{Jac-bou}, we know $h\leq \Jac(\mathcal {A})$, and thus the second part holds.
\qedd

\subsection{Algorithms to compute the differential Chow form}
Let $\mathcal {I}=\sat(\mathcal{A})$ be a prime differential ideal of dimension $d$  and $\mathcal {A}=\{A_{1},\ldots,A_{n-d}\}$ a given characteristic set of $\mathcal {I}$ under an arbitrary fixed ranking $\mathscr{R}$. In this section, we will give algorithms to compute the differential Chow form $F$ of $\mathcal {I}$ based on the order and degree bounds given in previous subsections.
Here, we use two different searching strategies by giving order and degree distinct priorities.

\subsubsection{Order priority}
In this section, we will give Algorithm 2 to compute the differential Chow form $F$ of $\mathcal {I}$
where the algorithm works adaptively by searching $F$ with order $h$ from $\ord(\mathcal{A})$ to $\Jac(\mathcal {A})$.
Indeed, by \cite{gao2010intersection}, we know that the order of $\mathcal {I}$ is equal to the maximum of all relative orders of $\mathcal {I}$,
and thus $\ord(F)\geq\ord(\mathcal{A})$, that is why  we start from $h=\ord(\mathcal {A})$.
 For a fixed order $h$, we search $F$ from $t=1$. If we cannot find $F$ with such a degree, then we repeat the procedure with $t+1$ until $t>\prod_{i=1}^{n-d}\deg(A_i)^{h+1}2^{(h+1)(d+1)}\big(2\sum_{i=1}^{n-d}(\deg(A_i)-1)+1\big)$.
 If for this $h$, $F$ cannot be found, then we repeat the procedure with $h+1$.
 In this way, we need only to handle problems with the real size and need not go to the upper bound in most cases.
Note that  the order bound given in Theorem~\ref{Jac-bou} and the degree bound given in Theorem~\ref{any-deg-bound} guarantee the termination of this algorithm.

\begin{algorithm}\label{alg-diffChowform2}
  \caption{\bf --- DChowform-2($\mathcal {A}$)} \smallskip
  \Inp{A characteristic set $\mathcal {A}=\{A_{1},\ldots,A_{n-d}\}$ of a prime differential ideal $\mathcal {I} $ under an arbitrary differential ranking $\mathscr{R}$}\\
  \Outp{The differential Chow form $F(\bu_0,\ldots,\bu_{d})$ of $\mathcal {I}$.}\medskip
  \noindent

  1. For $i=0,\ldots,d$, let
  $\P_i=u_{i0}+u_{i1}y_1+\cdots+u_{in}y_n$ and
  $\bu_i=(u_{i0},\ldots,u_{in})$.\\
  2. Set $h=\ord(\mathcal {A})$. \\
  3. Set $F=0$.\\
  4. While $F=0$ do\\
   \SPC 4.1. Set $t=1$, $\bv=\cup_{i=0}^{d}\bu_i^{[h]}$.\\
   \SPC 4.2. While $t\leqslant \prod_{i=1}^{n-d}\deg(A_i)^{h+1}2^{(h+1)(d+1)}\big(2\sum_{i=1}^{n-d}(\deg(A_i)-1)+1\big)$ do\\
   \SPC\SPC 4.2.1.  Set $F_0$ to be a homogenous GPol of degree $t$ in $\bv$.\\
   \SPC\SPC 4.2.2.  Set $\textbf{c}=\coeff(F_0,\bv)$. \\
   \SPC\SPC 4.2.3. Substitute $u_{i0}=-u_{i1}y_1-\cdots-u_{in}y_n\,(i=0,\ldots,d)$ into $F_{0}$ to get $F_1$.\\
  \SPC\SPC 4.2.4. Compute $F_2=\drem(F_1,\mathcal {A})$ under ranking $\mathscr{R}$.\\
   \SPC\SPC 4.2.5.  Set $\mathcal{P}=\coeff(F_2,\Theta(\Y)\cup\bv)$. Note $\mathcal{P}$ is a set of linear homogenous polynomials  
    \SPC \SPC\SPC\SPC in $\textbf{c}$.\\
   \SPC\SPC 4.2.6.  Solve the linear equation system $\mathcal{P}=0$.\\
   \SPC\SPC 4.2.7. If $\textbf{c}$ has a non-zero solution, then substitute it into $F_{0}$ to get $F$ and return $F$; 
            \SPC\SPC\SPC\SPC else $F=0$.\\
   \SPC\SPC  4.2.8.  $t:=t+1$.\\
   \SPC 4.3. h=h+1.\medskip
  \noindent

  /*/\; Pol and GPol stand for algebraic polynomial
  and generic  algebraic polynomial.\smallskip
   \noindent

   /*/\; $\coeff(F,V)$ returns the set of coefficients of $F$ as an algebraic  polynomial in
$V$.
\smallskip
\noindent

/*/\; $\drem(f, \mathcal{B})$ returns the differential remainder of $f$ w.r.t. an ascending chain $\mathcal{B}$.
\end{algorithm}

\begin{theorem}\label{com-2}
Let $\mathcal {I}=\sat(\mathcal {A})$ be a prime differential ideal of dimension $d$ and $\mathcal {A}=\{A_1,\ldots,\A_{n-d}\}$ a differential characteristic set under an arbitrary differential ranking. Set $m_i=\deg(A_i)$, $m=\max\{m_i\}$, $e_i=\ord(A_i)$, and $e=\max\{e_i\}$.
 Algorithm 2 computes the differential Chow form $F$ of $\mathcal {I}$ with at most
 $$O\big([n(m+1)]^{O(n(\Jac(\mathcal {A})+1)(2n-d+1)(e+d\Jac(\mathcal {A})+2\Jac(\mathcal {A})+d+1)}\big)$$ $\ff$-arithmetic operations.
\end{theorem}

\proof
By Theorem~\ref{Jac-bou} and  Theorem~\ref{any-deg-bound},
 Algorithm 2 computes a nonzero differential polynomial with minimal order and minimal degree under this order contained
in the differential ideal $[\sat(\mathcal{A}),\P_0,\ldots,\P_d]\cap\ff\{\bu_0,\ldots,\bu_d\}$, which is exactly the differential Chow form of $\sat(\mathcal{A})$.
So it remains to estimate the computational complexity of Algorithm 2.
Clearly, the complexity is dominated by steps 4.2.4 and 4.2.6.
Similarly as in the proof of Theorem \ref{com-1}, for fixed $h$ and $t$, step 4.2.4 and step 4.2.6 can be done with at most $O\big((n-d)(h+1)(m+1)^{O(n(n-d)(h+1)(e+h+1))}(2t+1)^{O(n(e+h+1))}\big)$ and  $O\big([(m+1)^{(n-d)(h+1)}(2t+1)]^{O((d+1)(n+1)(h+1)+n(e+h+1))}\big)$ arithmetic operations respectively.

  From Theorem~\ref{Jac-bou},  Step 4 may loop from $\ord(\mathcal {A})$ to $\Jac(\mathcal {A})$, and for each fixed $h$, step 4.2 may loop from 1 to $T(h)=\prod_{i=1}^{n-d}m_i^{h+1}2^{(h+1)(d+1)}\big(2\sum_{i=1}^{n-d}(m_i-1)+1\big)$. Thus, set $J=\Jac(\mathcal {A})$, the differential Chow form can be computed with at most
\begin{eqnarray}
& &\sum_{h=\ord(\mathcal {A})}^{J}\sum_{t=1}^{T(h)}\{O\big((n-d)(h+1)(m+1)^{O(n(n-d)(h+1)(e+h+1))}(2t+1)^{O(n(e+h+1))}+
\nonumber \\
& &[(m+1)^{(n-d)(h+1)}(2t+1)]^{O((d+1)(n+1)(h+1)+n(e+h+1))}\big)\}
\nonumber \\
&\leq&O\big(J\cdot T(J)\{(n-d)(J+1)(m+1)^{O((n-d)(J+1))}[(m+1)(2T(J)+1)]^{O(n(e+J+1))}
\nonumber \\
& &+[(m+1)^{(n-d)(J+1)}(2T(J)+1)]^{O(n(d+1)(J+1)+n(e+J+1))}\}\big)
\nonumber\\
&=&O\big(J\cdot T(J)[(m+1)^{(n-d)(J+1)}(2T(J)+1)]^{O(n(e+dJ+2J+d+1))}\big)
\nonumber
\end{eqnarray}
$\ff$-arithmetic operations. Here, to derive the above inequalities, $(m+1)^{(n-d)(J+1)}(2T(J)+1)$ $>n(e+dJ+2J+d+1)$ is assumed. Hence, the theorem follows by simply replacing $T(J)$ by the degree bound for $F$ given in Theorem ~\ref{any-deg-bound}.
\qedd

We use the following example to illustrate the above algorithm.
\begin{example}
Let $n=2$, $\mathcal {A}=\{y_2-y_1'\}$, $\mathscr{R}$ is the elimination ranking $y_1<y_2$. Clearly, $d=\dim(sat(\mathcal {A}))=1.$
We use this simple example to illustrate Algorithm 2. Let $\mathbb{P}_0=u_{00}+u_{01}y_1+u_{02}y_2$, $\mathbb{P}_1=u_{10}+u_{11}y_1+u_{12}y_2$, $\bu_0=(u_{00},u_{01},u_{02})$ and $\bu_1=(u_{10},u_{11},u_{12})$. In step 4.1, $h=\ord(\mathcal {A})=0$, $\textbf{v}=(u_{00},u_{01},u_{02},u_{10},u_{11},u_{12})$. In step 4.2, $t\leqslant 4$. We first execute steps 4.2.1 to 4.2.6 for $t=1$. Set $F_0=c_{01}u_{00}+c_{02}u_{01}+c_{03}u_{02}+c_{04}u_{10}+c_{05}u_{11}+c_{06}u_{12}$, $\textbf{c}=(c_{01}, \ldots, c_{06})$. In step 4.2.3, we get $F_1=-c_{01}u_{11}y_1-c_{01}u_{12}y_2+c_{02}u_{01}+c_{03}u_{02}-c_{04}u_{11}y_1-c_{04}u_{12}y_2+c_{05}u_{11}+c_{06}u_{12}$, and step 4.2.4 we get $F_2=-c_{01}u_{11}y_1-c_{01}u_{12}y_1'+c_{02}u_{01}+c_{03}u_{02}-c_{04}u_{11}y_1-c_{04}u_{11}y_1'+c_{05}u_{11}+c_{06}u_{12}$. Then $\mathcal {P}=0$ consists of equations $\{c_{01}=c_{02}=c_{03}=c_{04}=c_{05}=c_{06}=0\}$, $\mathcal {P}=0$ has a unique solution $\textbf{c}=(0,0,0,0,0,0)$. In step 4.2.8, $t=2$. Next we execute steps 4.2.1 to 4.2.6 for $t=2$. In the following computations, to save space, we will just list the number of equations and the solutions of the linear homogenous equation system $\mathcal {P}=0$ which are easily computed by Maple due to the strong sparsity of the system. For $t=2$, in step 4.2.5, we get 34 linear homogeneous polynomials in $\mathcal {P}$, and in step 4.2.6, we get $\mathcal {P}=0$ has a unique solution $\textbf{c}=(c_{01},\ldots,c_{21})=(0,\ldots,0)$. Next, we execute steps 4.2.1. to 4.2.6 for $t=3$. In step 4.2.5, we get 104 linear homogeneous polynomials in $\mathcal {P}$, and in step 4.2.6, we get $\mathcal {P}=0$ has a unique solution $\textbf{c}=(c_{01},\ldots,c_{56})=(0,\ldots,0)$. Then we execute steps 4.2.1 to 4.2.6 for $t=4$. In step 4.2.5, we get 259 linear homogeneous polynomials in $\mathcal {P}$, and in step 4.2.6, we get $\mathcal {P}=0$ has a unique solution $\textbf{c}=(c_{01},\ldots,c_{126})=(0,\ldots,0)$. Now,  in step 4.2, $t=5>4$. So we go on to Step 4.3 and obtain $h=1$.

 Since $F=0$, in step 4.1, set $\textbf{v}=(u_{00},u_{01},u_{02},u_{10},u_{11},u_{12},u'_{00},u'_{01},u'_{02},u'_{10},u'_{11},u'_{12})$ and $t=1$.
 Now, we execute Step 4.2 until $t> 8$ or $F\neq0$. We first execute steps 4.2.1 to 4.2.6 for $t=1$. $\mathcal {P}=0$ contains equations $\{c_{01}=c_{02}=c_{03}=c_{04}=c_{05}=c_{06}=c_{07}=c_{08}=c_{09}=c_{10}=c_{11}=c_{12}=0\}$, and has a unique solution $\textbf{c}=(c_{01},\ldots,c_{20})=(0,\ldots,0)$. Now, $t=2$ and we execute steps 4.2.1 to 4.2.6 for $t=2$. In step 4.2.5, we get 186 linear homogeneous polynomials in $\mathcal {P}$, and in step 4.2.6, $\mathcal {P}=0$ has a unique solution $\textbf{c}=(c_{01},\ldots,c_{78})=(0,\ldots,0)$. Then, we execute steps 4.2.1 to 4.2.6 for $t=3$. In step 4.2.5, we get 1122 linear homogeneous polynomials in $\mathcal {P}$, and in step 4.2.6, $\mathcal {P}=0$ has a unique solution $\textbf{c}=(c_{01},\ldots,c_{364})=(0,\ldots,0)$. Next, we execute steps 4.2.1 to 4.2.6 for $t=4$. In step 4.2.5, we get 5082 linear homogeneous polynomials in $\mathcal {P}$, and in step 4.2.6, $\mathcal {P}=0$ has a nonzero solution $\textbf{c}=(c_{01},\ldots,c_{1365})$ with $c_{110}=-q,c_{164}=q,c_{171}=q,c_{177}=-q,c_{256}=-q,c_{283}=q,c_{388}=q,c_{442}=-q,c_{449}=-q,c_{462}=q,c_{506}=q,
c_{568}=-q,c_{668}=q,c_{675}=-q,c_{725}=-q,
c_{760}=q$, where $q\in\mathbb{Q}$ and all the remaining $c$ equal to $0$.
Therefore, this algorithm returns $F=u_{00}u_{01}u_{11}u_{12}-u_{00}u_{02}u_{11}^2+u_{01}u_{02}u_{10}u_{11}-u_{01}^2u_{10}u_{12}+u_{00}'u_{02}u_{11}u_{12}-u_{00}'u_{01}u_{12}^2+u_{00}u_{02}u_{11}u_{12}'-u_{01}u_{02}u_{10}u_{12}'+u_{01}u_{02}u_{10}'u_{12}-u_{02}^2u_{10}'u_{11}
+u_{01}u_{02}'u_{10}u_{12}-u_{00}u_{02'}u_{11}u_{12}-u_{00}u_{02}u_{11}'u_{12}+u_{02}^2u_{10}u_{11}'-
u_{01}'u_{02}u_{10}u_{12}+u_{00}u_{01}'u_{12}^2$, which is the differential Chow form of $\mathcal {I}=\sat(\mathcal {A})$.
\end{example}

\subsubsection{Degree priority}
Algorithm 2 searches the differential Chow form with the order prior to the degree.
In other words, the output of Algorithm 2 is a nonzero  polynomial  in $[\sat(\mathcal{A}),\P_0,\ldots,\P_d]\cap\ff\{\bu_0,\ldots,\bu_d\}$ with minimal order and minimal degree under this order.
Thus, by the definition of differential Chow form, it must be the differential Chow form.

In this section, we give an alternative  algorithm to compute the differential Chow form of $\mathcal {I}=\sat(\mathcal{A})$ with the degree prior to the order during the searching strategy.
To be more precise, this algorithm works adaptively by searching $F$ from degree $t =1$ and for this fixed $t$ searching it with order $h$ from $\ord(\mathcal {A})$ to the order bound $\Jac(\mathcal{A})$. If a nonzero differential polynomial $F$ with degree $t$ is not found, then we repeat the procedures with degree $t+1$. If we find such an $F$, it requires to check whether $F$ is the differential Chow form.
We need to check it with the following conditions. Let $f \in \mathcal {F}\{\bu_0,\bu_1,\cdots,\bu_d\}$ be an irreducible differentially homogeneous polynomial, and $h=\ord(f)$. Let $\mathscr{R}_2$ be the elimination ranking $\cup_{i}\bu_i\backslash\{u_{i0}\}<\Y< u_{00}<\cdots<u_{d0} $ and $\mathscr{R}_2|_{\Y}=\mathscr{R}$.
Claim  (*):  $f$ is the differential Chow form if the following conditions are satisfied:
\begin{enumerate}[1.]
\item The differential remainder of $f$ w.r.t. $\{\mathcal {A}, \mathbb{P}_0,\ldots,\mathbb{P}_d\}$ under the ranking $\mathscr{R}_2$ is zero;
\item The differential remainder of each element in $\{\mathcal {A},\P_0,\ldots,\P_d\}$ w.r.t.  $\{f,\frac{\partial f}{\partial u_{00}^{(h)}}y_1-\frac{\partial f}{\partial u_{01}^{(h)}},\frac{\partial f}{\partial u_{00}^{(h)}}y_2-\frac{\partial f}{\partial u_{01}^{(h)}},\ldots,$
$\frac{\partial f}{\partial u_{00}^{(h)}}y_n-\frac{\partial f}{\partial u_{0n}^{(h)}}\}$ under the elimination ranking $\cup_{i}\bu_i< y_1<\cdots<y_n$   is zero;
while the differential remainder of $\init_{\mathcal {A}}\Sep_{\mathcal {A}}$ w.r.t. $\{f,\frac{\partial f}{\partial u_{00}^{(h)}}y_1-\frac{\partial f}{\partial u_{01}^{(h)}},\ldots,$
$\frac{\partial f}{\partial u_{00}^{(h)}}y_n-\frac{\partial f}{\partial u_{0n}^{(h)}}\}$ is nonzero.
\end{enumerate}
Let $\mathcal{J}=[\sat(\mathcal {A}), \mathbb{P}_0,\ldots,\mathbb{P}_d]\subset\ff\{\Y,\bu,u_{00},\ldots,u_{d0}\}$ where $\bu=\cup_{i}\bu_i\backslash\{u_{i0}\}$.
Before proving the claim, we first need the following two lemmas.

\begin{lemma}\label{check-1}
Let $f \in \mathcal {F}\{ \bu_0,\bu_1,\cdots,\bu_d\}$ be an irreducible differentially homogeneous polynomial, and $h=\ord(f)$.
Suppose the differential remainder of $f$ w.r.t. $\{\mathcal {A}, \mathbb{P}_0,\ldots,\mathbb{P}_d\}$ under the ranking $\mathscr{R}_2$ is zero. Set $\mathcal{C}=\{f,\frac{\partial f}{\partial u_{00}^{(h)}}y_1-\frac{\partial f}{\partial u_{01}^{(h)}},\ldots,\frac{\partial f}{\partial u_{00}^{(h)}}y_n-\frac{\partial f}{\partial u_{0n}^{(h)}}\}$,
then $\mathcal{C} \subseteq \mathcal {J}$.
\end{lemma}
\proof Obviously, $f \in \mathcal {J}$. Let $\xi =(\xi_1,\ldots,\xi_n)$ be a generic point of $\mathcal {I}=\sat(\mathcal{A})$ over $\mathcal {F}$ that is free from $\mathcal {F}\langle \bu\rangle $. Set $\eta_j=-\sum_{i=1}^n u_{ji} \xi_i$, then $(\xi_1,\ldots,\xi_n,\eta_0,\ldots,\eta_d)$ is a generic point of $\mathcal {J}$.
So $f(\bu,\eta_0,\ldots,\eta_d)=0$.   If we differentiate $f(\bu,\eta_0,\ldots,\eta_d)=0$ w.r.t. $u_{0\rho}^{(h)}(\rho=1,\ldots,n)$, then we have $\overline{\frac{\partial f}{\partial u_{0\rho}^{(h)}}}-\xi_{\rho}\overline{\frac{\partial f}{\partial u_{00}^{(h)}}}=0$, where $\overline{\frac{\partial f}{\partial u_{0\rho}^{(h)}}}$ and $\overline{\frac{\partial f}{\partial u_{00}^{(h)}}}$ are obtained by replacing $u_{00},\ldots,u_{d0}$ with $\eta_0,\ldots,\eta_d$ in $\frac{\partial f}{\partial u_{0\rho}^{(h)}}$ and $\frac{\partial f}{\partial u_{00}^{(h)}}$ respectively. So $\frac{\partial f}{\partial u_{00}^{(h)}}y_{\rho}-\frac{\partial f}{\partial u_{0\rho}^{(h)}} \in \mathcal {J}$ and therefore $\mathcal{C} \subseteq \mathcal {J}$.
\qedd

\begin{lemma}\label{check-2}
Let $f$ and $\mathcal{C}$ be as above. Set $\init_{\mathcal {A}}$ and $\Sep_{\mathcal {A}}$ to be the set of the initials and separants of $\mathcal {A}$ respectively. Suppose that the differential remainder of each element in $\{\mathcal {A},\P_0,\P_1,\ldots,\P_d\}$ w.r.t. $\mathcal{C}$ is zero, and the differential remainder of $\init_{\mathcal {A}}\Sep_{\mathcal {A}}$ w.r.t. $\mathcal{C}$ is nonzero. Then $\mathcal{C}$ is a characteristic set of $\mathcal {J}$ in $\mathcal {F}\{ \bu_0,\bu_1,\cdots,\bu_d,\mathbb{Y}\}$ w.r.t. the elimination ranking $\bu<u_{d0}<\ldots< u_{00}<y_1 < \ldots< y_n$.
\end{lemma}
\proof
Firstly, by the above lemma, $\mathcal{C} \subseteq \mathcal {J}$. Set $g_0=f$, $g_i=\frac{\partial f}{\partial u_{00}^{(h)}}y_{i}-\frac{\partial f}{\partial u_{0i}^{(h)}}(i=1,\ldots,n)$. Obviously, $\mathcal{C}$ is an irreducible auto-reduced set w.r.t. the elimination ranking $\bu<u_{d0}<\ldots< u_{00}<y_1 < \ldots< y_n$.
Thus, $\sat(\mathcal{C})$ is a prime differential ideal  with $\mathcal{C}$ being a characteristic set.
Therefore, it is sufficient to prove $\sat(\mathcal{C})=\mathcal {J}$.
 For any $g \in \mathcal {J}=[\sat(\mathcal {A}), \mathbb{P}_0,\ldots,\mathbb{P}_d)]$,
 we have $(\init_\mathcal {A}\Sep_\mathcal {A})^tg\in[\mathcal {A}, \mathbb{P}_0,\ldots,\mathbb{P}_d]$ for some $t \in \mathbb{N}$.
 Since the differential remainder of each $\P_i $ and each element in $\mathcal {A}$ w.r.t. $\mathcal{C}$ is zero,
 $\P_i\in\sat(\mathcal{C})$ and $\mathcal{A}\subset\sat(\mathcal{C})$.
 Thus, $(\init_\mathcal {A}\Sep_\mathcal {A})^tg \in \sat(\mathcal{C})$. Since $\sat(\mathcal{C})$ is a prime differential ideal, and the differential remainder of $\init_{\mathcal {A}}\Sep_{\mathcal {A}}$ w.r.t. $\mathcal{C}$ is nonzero, we have $g \in \sat(\mathcal{C})$ and it follows that $\mathcal {J}\subseteq \sat(\mathcal{C})$. Conversely, since $f$ is irreducible, we have $\frac{\partial f}{\partial u_{00}^{(h)}} \notin \mathcal {J}$. For, if not, then  $\frac{\partial f}{\partial u_{00}^{(h)}}\in\mathcal{J}\subseteq \sat(\mathcal{C})$,
and $\frac{\partial f}{\partial u_{00}^{(h)}}$ will be divisible by $f$, a contradiction.
So $\sat(\mathcal{C})=\big([\mathcal{C}]:(\frac{\partial f}{\partial u_{00}^{(h)}})^\infty\big)\subseteq \mathcal {J}$.
Thus, the lemmas is valid.
\qedd

\begin{pot}
By Lemmas~\ref{check-1} and \ref{check-2}, the claim is proved.
\end{pot}

With the above preparations, we now give Algorithm 3.
\begin{algorithm}[htdp]\label{alg-diffChowform3}
  \caption{\bf --- DChowform($\mathcal {A}$)} \smallskip
  \Inp{A characteristic set $\mathcal {A}=\{A_{1},\ldots,A_{n-d}\}$ of a prime differential ideal $\mathcal {I} $ under an arbitrary differential ranking $\mathscr{R}$}.\\
  \Outp{The differential Chow form $F(\bu_0,\ldots,\bu_{d})$ of $\mathcal {I}$.}\medskip
  \noindent

  1. For $i=0,\ldots,d$, let
  $\P_i=u_{i0}+u_{i1}y_1+\cdots+u_{in}y_n$ and
  $\bu_i=(u_{i0},\ldots,u_{in})$.\\
  2. Set $\widehat{h}=\Jac(\mathcal {A})$.\\
  3. Set $F=0$ and $t=1$.\\
  4. While $t\leqslant \prod_{i=1}^{n-d}\deg(A_i)^{\widehat{h} +1}2^{(\widehat{h}+1)(d+1)}\big(2\sum_{i=1}^{n-d}(\deg(A_i)-1)+1\big)$ do\\
  \SPC 4.1. Set $h=\ord(\mathcal {A})$.\\
   \SPC 4.2. While $h\leqslant \widehat{h}$ do\\
   \SPC\SPC 4.2.1.  Set $F_0$ to be a homogenous GPol of degree $t$ in $\bv=\cup_{i=0}^{d}\bu_i^{[h]}$.\\
   \SPC\SPC 4.2.2.    Set $\textbf{c}=\coeff(F_0,\bv)$. \\
   \SPC\SPC 4.2.3. Substitute $u_{i0}^{(k)}=-(u_{i1}y_1+\cdots+u_{in}y_n)^{(k)}\,(i=0,\ldots,d; 0\leqslant k\leqslant  h)$ into $F_{0}$ to  get $F_1$.\\
  \SPC\SPC 4.2.4. Compute $F_2=\drem(F_1,\mathcal {A})$ under ranking $\mathscr{R}$.\\
   \SPC\SPC 4.2.5.  Set $\mathcal{P}=\coeff(F_2,\Theta(\Y)\cup\bv)$. Note $\mathcal{P}$ is a set of linear homogenous polynomials      
   \SPC \SPC \SPC  in $\textbf{c}$.\\
   \SPC\SPC 4.2.6.  Solve the linear equation system $\mathcal{P}=0$.\\
   \SPC\SPC 4.2.7. If $\mathcal{P}=0$ has non-zero solutions, then pick one and substitute it  into $F_{0}$ to  get $F$;\\
    \SPC\SPC 4.2.8.  If $F\neq0$, then\\
    \SPC\SPC\SPC 4.2.8.1. If $F$ is not differentially homogeneous, then $F=0$, $\widehat{h}=h-1$, goto step  4.3.\\
   \SPC\SPC\SPC 4.2.8.2. For $1\leqslant i\leqslant n-d$, compute $\alpha _i=\drem(A_i,\mathcal{C}_F$), if $\alpha _i\neq 0$, then $F=0$, $\widehat{h}=h-1$,
   \SPC\SPC\SPC \SPC\SPC goto step 4.3, else $i=i+1$.\\
   \SPC\SPC\SPC 4.2.8.3. For $1\leqslant i\leqslant d$, compute $\beta _i=\drem(\mathbb{P}_i,\mathcal{C}_F)$, if $\beta_i\neq 0$, then $F=0$, $\widehat{h}=h-1$,  
   \SPC \SPC \SPC \SPC\SPC goto step 4.3, else $i=i+1$.\\
   \SPC\SPC\SPC 4.2.8.4.   Compute $\drem(\init_\mathcal {A}\Sep_\mathcal {A},\mathcal{C}_F)$, if it equals to zero, then $F=0$, $\widehat{h}=h-1$,  
   \SPC \SPC \SPC \SPC\SPC goto step 4.3.\\
   \SPC\SPC\SPC 4.2.8.5. Return $F$.\\
   \SPC\SPC  4.2.9.  $h:=h+1$.\\
   \SPC 4.3. t:=t+1.\medskip
  \noindent

 /*/\;$\mathcal{C}_F=\{F,\frac{\partial F}{\partial u_{00}^{(h)}}y_1-\frac{\partial F}{\partial u_{01}^{(h)}},\ldots,\frac{\partial F}{\partial u_{00}^{(h)}}y_n-\frac{\partial F}{\partial u_{0n}^{(h)}}\}$.\smallskip
   \noindent

  /*/\; Pol and GPol stand for algebraic polynomial
  and generic  algebraic polynomial.\smallskip
   \noindent

   /*/\; $\coeff(F,V)$ returns the set of coefficients of $F$ as an algebraic  polynomial in
$V$.\smallskip
\noindent

/*/\; $\drem(f, \mathcal{B})$ returns the differential remainder of $f$ w.r.t. an auto-reduced set $\mathcal{B}$.
\end{algorithm}

\begin{theorem}
Let $\mathcal {I}=\sat(\mathcal {A})$ be a prime differential ideal of dimension $d$ with $\mathcal {A}=A_1,\ldots,A_{n-d}$ a differential characteristic set under an arbitrary differential ranking. Set $m_i=\deg(A_i)$, $m=\max_i\{m_i\}$, $e_i=\ord(A_i)$, and $e=max_i\{e_i\}$.
 Algorithm 3 computes the differential Chow form of $\mathcal {I}$ with at most $$O\big([n(m+1)^{O((\Jac(\mathcal {A})+1)(2n-d+1))}]^{O(n(e+d\Jac(\mathcal {A})+2\Jac(\mathcal {A})+d+1))}\big)$$ $\ff$-arithmetic operations.
\end{theorem}
\proof  Firstly, we claim that ($\star$) for each fixed degree $t$,
step 4.2.9 will be executed if and only if $\mathcal{P}=0$ in 4.2.6 has only trivial solution $\textbf{0}$,
which implies that  for each fixed $t$, steps in 4.2.8 can be executed for at most one $h$.
Let $\mathcal{P}_{t,h}$ be the linear homogenous polynomial system obtained in step 4.2.5 for fixed degree $t$ and order $h$.
Suppose there exists an $h\leq \widehat{h}$ such that $\mathcal{P}_{t,h}=0$ has nonzero solutions
while $\mathcal{P}_{t,i}=0\,(i<h)$ has only zero solutions.
Take an arbitrary nonzero solution of $\mathcal{P}_{t,h}=0$ to obtain $F$.
If this $F$ does not satisfy  steps 4.2.8.1 to 4.2.8.4, then $F$ will be returned as an output.
Otherwise, $F$ is a nonzero differential polynomial in $\sat(\Chow(\mathcal{I}))= [\sat(\mathcal {A}),\mathbb{P}_0,\ldots,\mathbb{P}_d]\cap \mathcal {F}\{\bu_0,\ldots,\bu_d\}$
which is not the differential Chow form, so $\ord(\Chow(\mathcal{I}))\leq h$.
But if $\ord(\Chow(\mathcal{I}))= h$, $\Chow(\mathcal{I})$ divides $F$, a contradiction to  the fact that $\deg(\Chow(\mathcal{I}))>t.$
Thus, in this case, $\ord(\Chow(\mathcal{I}))<h$,  so just set $\widehat{h}=h-1$
and we do not need to execute step 4.2.9.

The algorithm aims to find a nonzero polynomial $F\in [\sat(\mathcal {A}),\mathbb{P}_0,\ldots,\mathbb{P}_d]\cap \mathcal {F}\{\bu_0,\ldots,\bu_d\}$
satisfying the conditions in Lemma~\ref{check-2} with minimal degree.
If such a polynomial  $F$ is found for a $(t,h)$, then it must be the differential Chow form.
Indeed, this $F$ must be irreducible, for $\mathcal{P}_{i,j}=0$ only possess zero solutions for $i<t$ and $j\leq h$.
By Lemma~\ref{check-2} and \cite[Lemma 4.10]{gao2010intersection}, $\mathcal{C}_f$ and $\mathcal{C}_{\Chow(\mathcal{I})}$ are both characteristic set of $\mathcal {J}= [\sat(\mathcal {A}),\mathbb{P}_0,\ldots,\mathbb{P}_d]\subset \mathcal {F}\{\Y,\bu_0,\ldots,\bu_d\}$ w.r.t. the elimination ranking $\bu<u_{d0}<\ldots< u_{00}<y_1< \ldots< y_n$, which implies that
 $F=a\cdot \Chow(\mathcal{I})$ for some $a\in\ff$ and so the output is just the differential Chow form of $\mathcal{I}$.

We now show such a polynomial can always be found.
In step 4.2.7, we just pick an arbitrary nonzero solution $\textbf{c}$ and substitute it into $F_{0}$ to get $F$, and then in step 4.2.8,  check if this $F$ satisfies the conditions described in Lemma~\ref{check-2}.
We claim that it is indeed enough to pick any one of the nonzero solutions in step 4.2.7.
Suppose there are two distinct solutions $\textbf{c}_1$ and $\textbf{c}_2$ of $\mathcal{P}=0$ obtained in step 4.2.6.
Let $F_1$ and $F_2$  be the polynomials obtained by substituting $\textbf{c}_1$ and $\textbf{c}_2$ into $F_0$ respectively.
Equivalently, we need to show that $F_1$ does not satisfy steps 4.2.8.1 to 4.2.8.4 if and only if  $F_2$ does not satisfy steps 4.2.8.1 to 4.2.8.4.
Suppose $F_1$ does not satisfy steps 4.2.8.1 to 4.2.8.4, then by  Lemma~\ref{check-2},  $F_1$ is the differential Chow form of $\sat(\mathcal{A})$.
Since $F_2$ has the same degree as $F_1$ and the same order guaranteed by claim ($\star$),
$F_2=a\cdot F_1\,(a\in\ff)$ must be the differential Chow form, which proves the claim.
By the above facts,  such a polynomial can always be found and  the output is the differential Chow form of $\mathcal{I}$.

We will estimate the complexity of the algorithm below.
 In each loop of step 4.2, the complexity of the algorithm is clearly dominated by step 4.2.4, step 4.2.6, and step 4.2.8.
 Similarly as in the proof of Theorem \ref{com-1}, for fixed $t$ and $h$, step 4.2.4. and step 4.2.6. can be done with at most $T_1=O\big((n-d)(h+1)(m+1)^{O(n(n-d)(h+1)(e+h+1))}(2t+1)^{O(n(e+h+1))}\big)$ and  $T_2=O\big([(m+1)^{(n-d)(h+1)}(2t+1)]^{O((d+1)(n+1)(h+1)+n(e+h+1))}\big)$ arithmetic operations respectively. In step 4.2.8.2, we need to compute the differential remainder of $A_i$ w.r.t. $\mathcal{C}_F$. By Lemma~\ref{alg-2}, this step can be done with at most $T_3=\sum_{i=1}^{n-d}(n+1)(e_i+1)(t+1)^{O(n(n+1)(e_i+1)(h+e_i+1))}(d_i+1)^{O(n(h+e_i+1))}$ arithmetic operations. Similarly, we get step 4.2.8.3 and step 4.2.8.4 can be done with at most $T_4=\sum_{i=1}^{n-d}(n+1)(e_i+1)(t+1)^{O(n(n+1)(e_i+1)(e_i+1))}3^{O(n(e_i+1))}$ and $T_5=2\sum_{i=1}^{n-d}(n+1)(e_i+1)(t+1)^{O(n(n+1)(e_i+1)(h+e_i+1))}(d_i+1)^{O(n(h+e_i+1))}$ arithmetic operations respectively. From Theorem ~\ref{any-deg-bound}, we know that step 4. may loop from 1 to $T=\prod_{i=1}^{n-d}m_i^{\Jac(\mathcal {A})+1}2^{(\Jac(\mathcal {A})+1)(d+1)}\big(2\sum_{i=1}^{n-d}(m_i-1)+1\big)$, and for each fixed $t$, from Theorem~\ref{Jac-bou},  Step 4.2 may loop from $\ord(\mathcal {A})$ to $\Jac(\mathcal {A})$. Thus, the differential Chow form can be computed with at most

\begin{eqnarray}
& &\sum_{t=1}^{T}\sum_{h=\ord(\mathcal {A})}^{\Jac(\mathcal {A})}(T_1+T_2+T_3+T_4+T_5)
\nonumber \\
&=&\sum_{t=1}^{T}\sum_{h=\ord(\mathcal {A})}^{\Jac(\mathcal {A})}O\big([(m+1)^{(n-d)(h+1)}(2t+1)]^{O((d+1)(n+1)(h+1)+n(e+h+1))}+
\nonumber \\
& &(n-d)(h+1)(m+1)^{O(n(n-d)(h+1)(e+h+1))}(2t+1)^{O(n(e+h+1))}+
\nonumber\\
& &3\sum_{i=1}^{n-d}(n+1)(e_i+1)(t+1)^{O(n(n+1)(e_i+1)(h+e_i+1))}(m_i+1)^{O(n(h+e_i+1))}+
\nonumber\\
& &\sum_{i=1}^{n-d}(n+1)(e_i+1)(t+1)^{O(n(n+1)(e_i+1)(e_i+1))}3^{O(n(e_i+1))}\big)
\nonumber\\
&\leq&O\big(T\cdot \Jac(\mathcal {A})\{[(m+1)^{(n-d)(\Jac(\mathcal {A})+1)}(2T+1)]^{O((d+1)(n+1)(\Jac(\mathcal {A})+1)+n(e+\Jac(\mathcal {A})+1))}+
\nonumber \\
& &(n-d)(n+1)(e+1)(T+1)^{O(n(n+1)(e+1)(e+\Jac(\mathcal {A})+1))}(d+1)^{O(n(\Jac(\mathcal {A})+e+1))}\}\big).
\nonumber
\end{eqnarray}
$\ff$-arithmetic operations. Here,  we assume $(m+1)^{(n-d)(\Jac(\mathcal {A})+1)}(2T+1)>n(e+d\Jac(\mathcal {A})+2\Jac(\mathcal {A})+d+1)$, and $\Jac(\mathcal {A})>>n$. Thus the theorem follows by simply replacing $T$ by the degree bound for $F$ given in Theorem ~\ref{any-deg-bound}.
\qedd

\begin{remark}
When using Algorithm 3 to compute the differential Chow form,
in step 4.2.8.1, we can examine whether the current nonzero differential polynomial $F$ satisfies the symmetric properties described in Theorem~\ref{th-chowproperty}.
If it is not symmetric, we can directly go onto Step 4.3.
\end{remark}

\begin{remark}
We use Figure 1 to illustrate the searching strategies of Algorithm 2 and Algorithm 3.
Both algorithms have their own advantages and defects in different situations.
Figure 1 shows Algorithm 2 has higher efficiency than Algorithm 3 in some cases.
And it may happen that Algorithm 3 has higher efficiency than Algorithm 2 in certain cases.
For example, let $n=2$ and $\mathcal {A}=\{(y'_1)^2y''_2-y_1\}$ with $\mathscr{R}$ being the elimination ranking $y_2<y_1$.
Here, the differential Chow form of $\sat(\mathcal{A})$ is of order $2$ and total degree $14$.
We use Figure 2 to show the steps which are needed to execute in Algorithm 2 and Algorithm 3 respectively for this example.
It is clear that  Algorithm 3 is of higher efficiency than Algorithm 2 in this particular example.

\begin{figure} \label{figure1}
     \begin{minipage}{0.5\textwidth}
    \begin{center}
          \begin{picture}(140,120)(0,0)
        \put(0,17){} \put(0,37){} \put(0,57){} \put(0,77){} \put(0,97){}\put(0,117){$t$}
        \put(16,0){} \put(36,0){} \put(56,00){} \put(76,0){} \put(96,0){}\put(116,0){$h$}
        \multiput(20,20)(20,0){6}{\line(0,1){100}}
        \multiput(20,20)(0,20){5}{$\uparrow$}
        \multiput(40,20)(0,20){3}{$\uparrow$}
        \multiput(20,20)(0,20){6}{\line(1,0){100}}
         \multiput(17.5,117.5)(0,-20){6}{$\circ $ }
         \multiput(37.5,57.5)(0,-20){3}{$\circ $ }
          \multiput(37.5,77.5)(0,-20){1}{$\bullet$ }
    \end{picture}
      \end{center}
\quad    {\qquad \qquad\qquad (a) Algorithm 2}
\label{fig-1}
  \end{minipage}
    \begin{minipage}{0.44\textwidth}
    \begin{center}
      \begin{picture}(140,120)(0,0)
        \put(0,17){} \put(0,37){} \put(0,57){} \put(0,77){} \put(0,97){}\put(0,117){$t$}
        \put(16,0){$$} \put(36,0){$$} \put(56,00){$$} \put(76,0){$$} \put(96,0){}\put(116,0){$h$}

        \multiput(20,20)(20,0){5}{$\rightarrow$}
                \multiput(20,40)(20,0){4}{$\rightarrow$}
                        \multiput(20,60)(20,0){3}{$\rightarrow$}
                                \multiput(20,80)(20,0){1}{$\rightarrow$}
         \multiput(20,20)(20,0){6}{\line(0,1){100}}
        \multiput(20,20)(0,20){6}{\line(1,0){100}}
         \multiput(17.5,77.5)(0,-20){4}{$\circ$}
           \multiput(37.5,57.5)(0,-20){3}{$\circ$  }
           \multiput(37.5,77.5)(0,-20){1}{$\bullet$ }
         \multiput(57.5,57.5)(0,-20){3}{$\circ$  }
         \multiput(77.5,57.5)(0,-20){3}{$\circ$ }
         \multiput(97.5,37.5)(0,-20){1}{$\star$  }
               \multiput(97.5,17.5)(0,-20){1}{$\circ$ }
         \multiput(117.5,17.5)(0,-20){1}{$\star$ }
      \end{picture}
    \end{center}
    \centerline{(b) Algorithm 3}
  \end{minipage}
  \caption{``$\circ$" means the algorithm is executed for the corresponding $(t,h)$ but $\mathcal{P}_{t,h}$ has only a zero solution,
   and ``$\star$" means  $\mathcal{P}_{t,h}$ has nonzero solutions but the corresponding nonzero $F$ is not the differential Chow form,
   while ``$\bullet$" means the corresponding $F$ is the output.}
  \end{figure}

\begin{figure} \label{figure2}
\subfigure{
\begin{minipage}[t]{0.5\textwidth}
   \centering
\begin{tikzpicture}
\coordinate (A) at (0.15,0.54);
\coordinate (B) at (0.15,2.5);
\coordinate (C) at (1.3,4.3);
\coordinate (D) at (2.3,1.4);
\coordinate (E) at (2.3,0.54);
\coordinate (J) at (0.15,0.6);
\fill[color=gray] (A)--(B)--(C)--(D)--(E)--cycle;
node[left] at (A) {(0,1)};
\node[left] at (B) {(0,60)};
\node[left] at (J) {(0,1)};
\node[right] at (C) {(1,720)};
\node[right] at (D) {(2,14)};
\node[right] at (E) {(2,1)};
\tikz\draw[thick,<->,>=stealth](0,4)node[above=3pt] {$t$}-|(0,0)--(4,0)node[right] {$h$};

\end{tikzpicture}
 \centerline{(c) Algorithm 2}
 \end{minipage}}
\subfigure{
\begin{minipage}[t]{0.5\textwidth}
\centering
\begin{tikzpicture}
\coordinate (F) at (0.15,0.54);
\coordinate (G) at (0.15,1.4);
\coordinate (H) at (2.3,1.4);
\coordinate (I) at (2.3,0.54);
\fill[color=gray] (F)--(G)--(H)--(I)--cycle;

\node[left] at (F) {(0,1)};
\node[left] at (G) {(0,14)};
\node[right] at (H) {(2,14)};
\node[right] at (I) {(2,1)};
\tikz\draw[thick,<->,>=stealth](0,4)node[above=3pt] {$t$}-|(0,0)--(4,0)node[right] {$h$};
\end{tikzpicture}
 \centerline{(d) Algorithm 3}
 \end{minipage}
 }
 \caption{Both Algorithm 2 and Algorithm 3 return a differential polynomial $F$ with $h=2$ and $t=14$.
 Algorithm 2 is executed at all the integer lattice points $(h,t)$ which lie in the gray convex polygon as shown in the figure (c),
  while Algorithm 3 is executed at all the integer lattice points $(h,t)$ of the gray convex polygon in the figure (d).}
 \end{figure}
\end{remark}

We conclude this section by giving an application of the algorithms in this paper.
Given a characteristic set $\mathcal{A}$ of a  prime differential ideal $\mathcal {I}$ under an arbitrary ranking,  Theorem~\ref{Jac-bou} shows that $\ord(\mathcal{I})\leq \Jac(\mathcal{A})$. But what is the precise order of $\mathcal{I}$? And how to compute it?

Since $\ord(\CI)=\ord(\Chow(\CI))$, if the differential Chow form of $\mathcal{I}$ has been computed, then clearly we can read off the order of $\mathcal{I}$. Thus, the above problem can be solved by computing the differential Chow form of $\CI$.
\section{Conclusion}

In this paper, we propose algorithms to compute differential Chow forms for prime differential ideals represented by their characteristic sets under arbitrary rankings
and estimate the computing complexities of these algorithms.
In general, two cases are considered according to whether the given ranking is orderly or not.

For a prime differential ideal  given by its characteristic set  under some orderly ranking, we first estimate the degree bounds for  the differential Chow form and
then based on the degree bound and also the precise order, we compute the differential Chow form  with linear algebraic techniques.
For a prime ideal given by characteristic sets under arbitrary rankings, we first give the order bound for the differential Chow form which is the Jacobi number of the
characteristic set. Then with the degree bound similar to that in the first case, we devise an algorithm to compute its differential Chow form.
Both algorithms need single exponential number of arithmetics in the worst case.

Recent study in differential algebra owes to the idea of using a wider class of differential ideals than prime differential ideals, namely, characterizable differential ideals \cite{hubert2000factorization}.
It is interesting to  compute differential Chow forms for characterizable differential ideals.
For, once the differential Chow form has been computed, by factoring the differential Chow form, we can give the irreduandant irreducible decomposition for the original differential ideal. The main difficulty in the generalization is to extend Theorem~\ref{Jac-bou} from differential prime ideals to differential characterizable ideals
under non-orderly rankings.

\section*{References}
 \bibliographystyle{model1b-num-names}
\bibliography{reference}

\end{document}